\RequirePackage{fix-cm}

\documentclass[leqno, fleqn, centertags, 12pt]{article}

\usepackage{latexsym}
\usepackage{amsmath}
\usepackage{amscd}
\usepackage{amssymb}
\usepackage{verbatim}

\usepackage[T1]{fontenc}
\usepackage{latexsym}
\usepackage{manfnt}
\usepackage[sans]{dsfont}
\usepackage{palatino}
\usepackage[sc, osf]{mathpazo}
\usepackage{eulervm}
\usepackage{euscript}

\usepackage{textcomp}

\usepackage{fullpage}

\newskip\nineskipamount \nineskipamount=9pt plus 0pt minus 0pt
\newskip\zeroskipamount \zeroskipamount=0pt plus 0pt minus 0pt
\usepackage[noorphans,vskip=\nineskipamount]{quoting}

\makeatletter
\renewcommand{\@makefntext}[1]{\vspace*{0.5ex}\parindent=0em
\hspace*{-0.4em}
\hbox to 0.4em{\hss\@makefnmark}\hspace*{0.4em}{#1}
}
\makeatother

\newcounter{mysectionnumber}
\newcommand{\mysection}[2]{
\refstepcounter{mysectionnumber}
\section*{ \textnormal{{\themysectionnumber.} {#1}}}\label{#2}}

\newcommand{\mynonumbersection}[2]{
\vspace{-0.0ex}
\section*{{}\hspace*{0.00em}$\phantom{1.}$\textnormal{{#1}}}\label{#2}\vspace{\bigskipamount}}  

\newcommand{\myit}[1]{\textbf{\textit{#1}}\hspace{0.0em}}

\newcounter{myparnum}[mysectionnumber]
\renewcommand{\themyparnum}{\themysectionnumber.\arabic{myparnum}}
\newcommand{\mypar}[2]{\refstepcounter{myparnum}{\vspace{\medskipamount}\textbf{{\themyparnum. #1}\label{#2}}\hspace{0.5em}}}

\newcommand{\myuppar}[1]{\vspace{\medskipamount}\textbf{#1}\hspace*{0.5em}}

\newcommand{\proof}{\vspace{\medskipamount}{\textbf{{\emph{Proof}.}}\hspace*{0.7em}}}
\newcommand{\eproof}{ $\blacksquare$}

\newcommand{\dis}{\displaystyle}

\def\sss{\hspace{0.05em}\ }
\def\dss{\hspace{0.1em}\ }
\def\trs{\hspace{0.15em}\ }
\def\qss{\hspace{0.2em}\ }
\def\pss{\hspace{0.3em}\ }
\def\oss{\hspace{0.4em}\ }

\def\halfff{\hspace*{0.025em}}
\def\fff{\hspace*{0.05em}}
\def\dff{\hspace*{0.1em}}
\def\trf{\hspace*{0.15em}}
\def\qff{\hspace*{0.2em}}
\def\pff{\hspace*{0.3em}}
\def\off{\hspace*{0.4em}}
\def\pquad{\hspace*{1.5em}}

\def\ttff{{\hspace*{-0.05em}--\hspace*{0.15em}}}

\newcommand{\hnsp}{\hspace*{-0.05em}}
\newcommand{\nsp}{\hspace*{-0.1em}}
\newcommand{\nnsp}{\hspace*{-0.15em}}

\newcommand{\dnsp}{\hspace*{-0.2em}}

\renewcommand{\geq}{\geqslant}
\newcommand{\id}{\mathop{\mbox{id}}\nolimits}

\newcommand{\cc}{\mathbold{C}}
\newcommand{\ccc}{\mathbb{C}}
\newcommand{\rr}{\mathbold{R}}
\newcommand{\rrr}{\mathbb{R}}
\newcommand{\rtwo}{\mathbold{R}^{2}}
\newcommand{\uu}{\mathbb{U}}

\newcommand{\re}{\mbox{Re}\qff}
\newcommand{\im}{\mbox{Im}\qff}

\newcommand{\toto}{\longrightarrow}
\newcommand{\ttoo}{\hspace*{0.2em}\longrightarrow\hspace*{0.2em}}

\begin{document}

\setlength{\baselineskip}{12pt plus 0pt minus 0pt}
\setlength{\parskip}{12pt plus 0pt minus 0pt}
\setlength{\abovedisplayskip}{12pt plus 0pt minus 0pt}
\setlength{\belowdisplayskip}{12pt plus 0pt minus 0pt}

\newskip\smallskipamount \smallskipamount=3pt plus 0pt minus 0pt
\newskip\medskipamount   \medskipamount  =6pt plus 0pt minus 0pt
\newskip\bigskipamount   \bigskipamount =12pt plus 0pt minus 0pt

\title{\vspace*{-0ex}{The geometric meaning of the complex dilatation}}
\date{}
\author{Nikolai V. Ivanov}

\footnotetext{\hspace*{-0.65em}\copyright\ 
Nikolai V. Ivanov,\oss 2017,\qss 2023.\qss}

\maketitle

\vspace*{3pt}\vspace*{12pt}
\myit{\hspace*{0em}\large Contents}\vspace*{1ex} \vspace*{\bigskipamount}\\ 
\myit{Introduction} \hspace*{0.5em} \vspace*{1ex}\\
\myit{1.}\hspace*{0.5em} Quadratic forms and conformal structures on $\cc$ \hspace*{0.5em} \vspace*{0.25ex}\\
\myit{2.}\hspace*{0.5em} The Poincar\'e invariant of a quadratic form \hspace*{0.5em} \vspace*{0.25ex}\\
\myit{3.}\hspace*{0.5em} The dilatation of real linear maps\qss $\cc\toto\cc$\hspace*{0.5em}\vspace*{0.25ex}\\
\myit{4.}\hspace*{0.5em} The action of real linear maps\qss $\cc\toto\cc$\hspace*{0.5em}\vspace*{0.25ex}\\
\myit{5.}\hspace*{0.5em} Complex vector spaces of dimension $1$  \hspace*{0.5em} \vspace*{1ex}\\
\myit{References}\hspace*{0.5em}  \hspace*{0.5em}  \vspace*{0.25ex}

\renewcommand{\baselinestretch}{1}
\selectfont

\mynonumbersection{Introduction}{introduction}

\myuppar{The classical complex dilatation.} 
As usual,\qss we will identify $\cc$ with $\rtwo$ by the map 
$x\qff +\qff i\dff y \longmapsto (x\fff,\pff y) $\dnsp.\oss 
Let $U$ be an open subset of $\cc$ and let $f\colon U \toto \cc$ be a map 
differentiable as a map from $U\dff \subset\dff \cc\qff =\qff \rtwo$ to $\rtwo$\dnsp.\oss 
The complex dilatation $\mu_{\dff f}$ of $f$ is a measure
of the distortion of the conformal structure of
$\cc$ by $f$\dnsp.\oss 
Let us recall its definition following Ahlfors\qss \cite{a}.\qff\oss 
Let\sss $z\qff =\qff x\qff +\qff i\dff y$\sss and let\qss 
\[
\quad
f\fff(z)\off =\off u\fff(x\fff,\pff y)\qff +\qff i\dff v\fff(x\fff,\pff y)\dff.
\] 
Let\oss
$\dis 
f_x\off =\off u_{\dff x}\qff  +\qff i\dff v_{\dff x}$\oss and\oss
$\dis
f_y\off =\off u_{\dff y}\qff  +\qff i\dff v_{\dff y}$\oss
be the usual partial derivatives and let
\[
\quad
f_{\fff z}\off =\off 
\bigl(f_x\qff -\qff i\dff f_y\bigr)\bigl/2,\quad\quad 
f_{\dff\overline{z}}\off =\off 
\bigl(f_x\qff +\qff i\dff f_y\bigr)\bigl/2\fff.
\]
The\qss \emph{complex dilatation}\qss of\dss $f$\sss is defined as the quotient\sss 
$\mu_{\dff f}\off =\off f_{\dff\overline{z}}\left/\dff f_{\fff z}\right.$\nsp.\qss

It is assumed that $f$ has non-zero Jacobian and is orientation-preserving. 
This ensures that\sss $f_z\qff \neq\qff 0$\sss 
and,\qss moreover\halfff,\qss that\sss $|\dff f_z\dff|\qff >\qff |\dff f_{\overline{z}}\dff|$\sss 
and hence\sss $|\dff\mu_{\dff f}\dff|\qff  <\qff 1$\nnsp.\oss 
In other words,\dss $\mu_{\dff f}$ belongs to the open unit disc $\mathbb{U}$\dss
in the complex plane.\oss
See\qss \cite{a}\fff,\oss Section\qss 1.A.\qss
We will sometimes refer to $\mu_{\dff f}$
defined as above as the\qss \emph{classical complex dilatation}.

Ahlfors explains the geometric meaning of $\mu_{\dff f}$ as follows. 
At any given point of $U$ the tangent map\sss $T_{\dff f}$\sss 
of\sss $f$\sss is a\qss \emph{real}\qss linear map.\oss
Since $f$ has non-zero Jacobian,\qss $T_{\dff f}$ is invertible, and therefore 
takes circles with the center at the origin into an ellipses.\oss 
The ratio of the major axis to the minor axis of every such ellipse 
is equal to the\qss \emph{dilatation}\oss
\[
\quad
D_{\dff f}
\off =\off
\frac{|\dff f_{\fff z}\dff |\qff  +\qff  |\dff f_{\dff\overline{z}}\dff |}
{\qff |\dff f_{\fff z}\dff |\qff  -\qff  |\dff f_{\dff\overline{z}} \dff|\qff}\qff.
\]
of\sss $f$\nnsp.\oss
It is related to the absolute value\sss
$d_{\dff f} \qff =\fff\qff |\dff \mu _{\dff f}\dff |$\sss
of\sss $\mu _{\dff f}$\sss by the formulas
\begin{equation}
\label{dilatation}
\quad
D_{\dff f}
\off =\off 
\frac{1\qff +\qff d_{\dff f}}{\dff 1\qff -\qff d_{\dff f} \dff}\qff,
\hspace*{2em} 
d_{\dff f}
\off =\off 
\frac{D_{\dff f}\qff -\qff 1}{\dff D_{\dff f}\qff +\qff 1 \dff}\qff.
\end{equation}
The ratio\sss
$|\dff T_{\dff f}\dff(z)\dff |\qff\bigl/\qff|\dff z\dff |\bigr.
$\sss
is maximal when\sss
$\arg z\off =\off
\arg \mu_{\dff f}\bigl/2$\nsp.\oss

Therefore,\oss the complex dilation\dss $\mu_{\dff f}$\dss is characterized 
by the following two properties.\vspace{\medskipamount}

(i)\quad\hspace{0em} \emph{$d_{\dff f}\off =\off 
|\dff \mu_{\dff f}\dff |$\dss
is related to the ratio\sss $D_{\dff f}$\sss of\dss the axes of an image ellipse
by\qss \textup{(\ref{dilatation}).}}
 
(ii)\quad\fff \emph{$\arg \mu_{\dff f}\off =\off 2\fff\alpha$\nnsp,\qss
where\dss $\alpha$\sss is the direction of the maximal distortion of\sss $T_{\dff f}$\nsp.} 
\vspace{\medskipamount}

This characterization is usually considered as 
the explanation of the geometric meaning of\dss $\mu_{\dff f}$\nnsp.\oss
Indeed\halfff,\qss 
\emph{the\dss absolute\dss value\sss 
$|\dff \mu_{\dff f}\dff |$\sss 
encodes\dss the\dss shape}\qss 
of the ellipses to which $T_{\dff f}$\sss 
takes the circles with the center at the origin.\oss
These ellipses are circles if and only if\dss 
$\mu_{\dff f}\qff =\qff 0$\dnsp,\qss
and if\dss 
$\mu_{\dff f}\qff \neq\qff 0$\dnsp,\qss 
then\qss 
\emph{the\dss argument\dss $\arg \mu_{\dff f}$\sss
together\dss with}\dss $T_{\dff f}$\sss 
determines the\qss \emph{directions}\qss
of the axes of such an ellipse.\oss
Namely,\oss the direction of the major axis 
is equal to the image under $T_{\dff f}$ of the direction
with the angle $(\arg \mu_{\dff f})/2$\nnsp.\oss

Still\halfff,\qss the most\sss natural way to encode the shape of an ellipse
is to take the ratio of its axes,\oss and\sss $D_{\dff f}$\sss 
seems to be
a much more natural measure of the distortion of shapes by $f$ than
the absolute value of the complex dilatation.\qss
The replacement of\dss $D_{\dff f}$ by\sss $d_{\dff f}$\sss 
is not justified geometrically at all\halfff.\qss

The main goal of this paper is to propose a geometric framework in which the complex dilatation
$\mu_{\dff f}$ and\halfff,\qss in particular\halfff,\qss 
its absolute value $d_{\dff f}$ appear naturally.\qss
This requires a radical change of the point of view:\qss 
instead of considering the individual ellipses,\qss
one should consider the space of all ellipses with the center at the origin
up to homothetic transformations centered at the origin.\qss
In other terms,\qss one needs to consider the\qss \emph{space of all conformal structures on}\dss $\cc$\dss
(considered as a tangent space to itself\dff).\qss

\myuppar{The argument of the complex dilatation.}
The angle $\alpha$ is defined 
only up to addition of\sss $\pm\qff\pi$\nnsp.\qss 
Replacing it by $2\fff\alpha$ is a way to kill this indeterminacy.\qss
We will see that the expression
$2\fff\alpha$ also appears naturally.\oss 

But\sss $\mu_{\dff f}$\sss 
determines the direction of the image ellipse
only together with $T_{\dff f}$\nnsp.\oss
In fact\halfff,\qss the\qss \emph{preimages}\qss 
of circles with the center $0$ under\sss $T_{\dff f}$\sss
are also ellipses and\sss $\alpha$\sss 
is the direction of the minor axes of these ellipses.\oss
Therefore,\qss $\mu_{\dff f}$\sss determines the direction
of the preimage,\qss rather than image,\qss ellipses.\oss 
This crucial feature is rarely mentioned.\qss

\myuppar{An approach to the geometric interpretation of the complex dilatation.}
By the definition\dss $\mu_{\dff f}$\dss at a point depend
only on the tangent map\dss $T_{\dff f}$\dss at this point.\oss 
In fact\halfff,\qss one can define the complex dilatation\dss $\mu_{\dff T}$\dss of any\dss
\emph{orientation-preserving\dss real\dss linear}\qss map\qss $T\colon\cc\toto\cc$\qss in such a way that\qss 
$\mu_{\dff f}\off =\off \mu_{\dff T_{\dff f}}$\qss for any differentiable map $f$ as above.\oss
Namely,\qss let\vspace{-1.5pt}
\[
\quad
\left(
\begin{array}{cc} u_{\dff 1} & u_{\dff 2} \\ 
v_{\fff 1} & v_{\fff 2} 
\end{array} 
\right)
\]

\vspace{-12pt}\vspace{-1.5pt}
be the matrix of $T$ with respect to the standard basis of\qss 
$\cc\qff =\qff \rtwo$\dnsp.\oss
Let\vspace{0pt}
\[
\quad
T_{\dff 1}
\off  =\off  
u_{\dff 1}\qff +\qff i\dff  v_{\fff 1}\dff,
\quad\quad 
T_{\dff 2}
\off  =\off  
u_{\dff 2}\qff +\qff i\dff  v_{\fff 2}\dff,
\hspace*{1.5em}\mbox{and}
\]

\vspace{-37.5pt}
\[
\quad
T_{\dff z}
\off  =\off  
\bigl(\dff T_{\dff 1}\qff -\qff i\dff  T_{\dff 2}\fff\bigr)\bigl/2\dff,
\quad\quad
T_{\dff \overline{z}}
\off  =\off  
\bigl(\dff T_{\dff 1}\qff +\qff i\dff  T_{\dff 2}\fff\bigr)\bigl/2\dff ,
\]

\vspace{-12pt}\vspace{0pt}
and define the\qss \emph{complex\dss dilatation}\qss 
of\sss $T$ as\sss 
$\mu_{\dff T}
\off =\off 
T_{\dff \overline{z}}\qff\bigl/\qff T_{\dff z}\bigr.$\nnsp.\qss
Clearly,\dss 
$\mu_{\dff f}\off =\off \mu_{\trf T_{\dff f}}$ for any differentiable map $f$ as above,\qss 
and Alhfors's interpretation
of $\mu_{\dff f}$ applies mutatis mutandis to $\mu_{\dff T}$\nnsp.\oss
Therefore,\qss the real linear maps $\cc\toto\cc$
are at heart of the matter\halfff.\oss 
Our approach is based on the following observations.\vspace{-1.5pt}

\myuppar{\dnsp$\bullet$\dnsp} 
A natural measure of distortion of the conformal structure 
of\dss $\cc$\dss by a\qss \textit{real\dss linear\dss map}\qss 
$T\colon\cc\toto\cc$\qss 
is the pull-back by $T$ of the standard conformal structure on\dss $\cc$\nnsp.\qss\vspace{-1.5pt}

\myuppar{\dnsp$\bullet$\dnsp}
The set $U\fff(\cc)$ of all conformal structures on\ $\cc$\dss
carries a canonical structure of a model of the hyperbolic plane.\oss
In fact\halfff,\dss $U\fff(\cc)$ can be naturally identified with the open unit disc $\mathbb{U}$
together with its structure of\sss the Klein model of the hyperbolic plane.\oss\vspace{-1.5pt}

\myuppar{\dnsp$\bullet$\dnsp} 
The standard isomorphism of the Klein model 
with the Poincar\'e unit disc model transforms this measure of distortion
into the classical complex dilatation.
\vspace{6pt}\vspace{-1.5pt}

The angle $2\fff\alpha$ naturally appears already in the Klein model,\oss
and the replacement of\dss $D_{\dff f}$\qss by\sss $d_{\dff f}$ results from
the transition from the Klein to the Poincar\'e model.\oss

\myuppar{Complex dilatation and\dss Arnold's\dss approach\sss to altitudes
in\sss hyperbolic\sss geometry.}
A\qss \emph{conformal\sss structure}\qss on $\cc$\sss is\dss defined as
a\qss (positive or negative)\qss definite quadratic form on $\cc$
considered up\sss to a non-zero real\sss factor\halfff.\qss
So,\qss on of\dss the key\sss elements of\sss our approach\sss is\dss
the canonical\sss structure of\sss the Klein model of the hyperbolic plane
on\sss the space of\dss definite quadratic forms on $\cc$ considered
up\sss to a non-zero real\sss factor\halfff.\qss
This\sss is\sss also a key element\sss of\dss Arnold's\dss approach\sss
to\sss the\sss triangle altitudes\sss theorem 
in\dss Lobachevsky\qss (hyperbolic)\qss geometry\qss \cite{ar}.\qss
The author has\sss to admit\sss that\sss he missed\sss this fact\sss
even after writing an exposition of\dss Arnold's\dss ideas\qss \cite{i}\qss
and completing\sss the first\sss version of\dss the present\sss paper\halfff.\qss
({\fff}In\qss \cite{i}\qss the\dss Poisson\dss bracket\sss of\dss quadratic forms on
$\cc\qff =\qff \rtwo$
was replaced\sss by\sss the commutator of\dss 
traceless $2\times 2$ real\sss matrices.)\vspace{-0.375pt}

\myuppar{The action of real linear maps.}
Let\qss $u\dff\colon\dff \cc\toto \cc$\qss is a real linear orientation-preserving isomorphism.\oss
Then there is a unique map\qss $u^*\dff\colon\dff \mathbb{U}\toto \mathbb{U}$\qss
such that\sss
$\mu_{\dff T\dff\circ\dff u}
\qff =\qff
u^*\dff(\mu_{\dff T})$\sss
for any orientation-preserving real linear map\qss $T\colon\cc\toto\cc$\nnsp.\oss
Moreover\halfff,\qss $u^*\dff\colon\dff \mathbb{U}\toto \mathbb{U}$\qss
is complex-analytic.\oss
Let us quote\qss J.H.\qss Hubbard\qss \cite{h},\oss p.\oss 162:\vspace{-9pt}

\begin{quoting}
The fact that\sss $u^*$\sss is analytic has far-reaching consequences.\qss
When you dig down to where the complex analytic structure of\dss 
Teichm\"{u}ller space comes down from\qss
(a highly nontrivial result\halfff,\qss
with a rich and contentious history\halfff,\oss
involving Ahlfors,\oss Rauch,\oss Grothendieck,\oss Bers,\oss and many others),\oss
you will find that this is the foundation of it all\qss 
(see\sss \ldots).\oss
Thus do little acorns into mighty oak trees grow.                                                                                                                                                                                                                                                                                                                                                                                                                                                                                                                                                                                                                                                                                                                                                                                                                                                                                                                            
\end{quoting}

\vspace{-9pt}
Uniqueness of\dss $u^*$\dss is straightforward.\oss
The complex-analyticity can be proved by a direct calculation.\qss
See\qss \cite{h},\oss Proposition\qss 4.8.10.
The approach of this paper offers a conceptual explanation of the complex-analyticity\halfff.qss
Taking the pull-backs of conformal structures by $u$\sss leads
a map\sss $U\fff(\cc)\toto U\fff(\cc)$\dnsp.\qss
By the very definition of the Klein model of the hyperbolic plane,\qss
this induced map is its automorphism\qss
(i.e.\qss preserves\sss lines,\qss angles,\qss distances,\qss etc.).\qss
Passing to the Poincar\'e model turns this map
into an automorphism\sss
$u^*\dff\colon\dff \mathbb{U}\toto \mathbb{U}$\sss 
of\dss the latter\halfff.\qss
The property\sss 
$\mu_{\dff T\dff\circ\dff u}
\qff =\qff
u^*\dff(\mu_{\dff T})$\sss is immediate,\qss
and hence this recovers the classical map $u^*$\nnsp.\qss
On the other hand,\oss it is well known that the automorphisms
of the Poincar\'e model are complex-analytic.\qss
See\dss Section\qss \ref{action}\qss for\sss the details.\vspace{-0.375pt}

\myuppar{The paper\halfff.}
The main part of the paper\halfff,\qss 
namely Sections\qss \ref{quadratic forms}\off--\off\ref{action},\qss
is devoted to a theory of the complex dilatation of real linear maps\qss
$\cc\toto\cc$\qss based on the above observations.\oss
The main result\dss is\dss Theorem\qss \ref{p-equal-classical}.\qss
In Section\qss \ref{one-dim-over-cc}\qss 
this theory is used to deal with 
real linear maps\qss $V\toto W$\qss
from one complex vector space of dimension $1$ to another\halfff.\qss
The complex dilatation of such a map is not a number\halfff,\qss
but a tensor\halfff.\oss
The author resisted the temptation to present this
theory in such generality from the very beginning.\oss 

The present paper grow out of the observation that the replacement 
of\sss $D_{\dff f}$ by $d_{\dff f}$ corresponds to the
transition from the Klein to the Poincar\'e model of the hyperbolic plane.\oss
An preliminary exposition\qss \cite{i95}\qss of this idea was written back in 1995\dss
in response to a question of\qss J.D.\qss McCarthy  
about the geometric meaning of the complex dilatation.\qss
I am grateful to him for asking\sss this question,\qss his interest,\qss and many conversations.

\newpage
\mysection{Quadratic\qss forms\qss and\qss conformal\qss structures\qss on\qss $\cc$}{quadratic forms}

\vspace{6pt}
\myuppar{Conformal structures.} 
Recall\sss that\sss a\dss \emph{quadratic form on a real vector space}\dss $V$ 
is a function $q\colon V\toto\rr$ such that\dss
$q\fff(v)\qff =\qff B\fff(v\fff,\pff v) $\dss for some\dss 
\emph{symmetric bilinear form}\dss $B\colon V\times V\to\rr$\dss and all\dss $v\in V$\dnsp.\oss 
The bilinear form $B$ can be reconstructed from the quadratic form $q$ as
its\qss \emph{polarization}\dss 
$B_q\fff(v\fff,\pff w)
\off =\off 
(q\fff(v\qff +\qff w)\qff -\qff q\fff(v)\qff -\qff q\fff(w)\fff)/2$\nnsp.\oss
A quadratic form $q$ is\qss \emph{positive definite}\qss 
if $q(v)>0$ for all $v\qff \neq\qff 0 $\dnsp,\oss
and\qss \emph{negative definite}\qss if $q(v)<0$ for all $v\qff \neq 0\qff $\dnsp.\oss 
It is\qss \emph{definite}\qss if\dss it\dss is\dss
either positive or negative definite.\qss
A\dss 
\emph{quadratic form on a complex vector space}\dss $V$ is\dss a quadratic form
on $V$ considered as a\qss \emph{real vector space}.\oss

A\dss \emph{conformal structure}\dss on a real or complex vector space $V$ 
is a definite quadratic form on $V$
considered up to multiplication by a non-zero real number\halfff,\oss or\halfff,\oss
equivalently, as a positive definite quadratic
form considered up to multiplication by a positive real number.
The conformal structure determined by a 
quadratic form $q$ is called the\dss \emph{conformal class of} $q$
and denoted by $[\fff q\dff]$\nnsp.\qss
The set $U\fff(V)$ of all conformal structures on $V$ is a subspace of the\qss
\emph{projective space}\dss $\mathbb{P}\hnsp Q(V)$ associated with the\qss 
vector\qss \emph{space\dss $Q(V)$\sss of\dss quadratic\sss forms on}\dss $V $\dnsp.\oss 
Namely,\qss $U\fff(V)$\dss is the image in $\mathbb{P}\hnsp Q(V)$ of the set of all definite forms.\oss

We are interested only in\sss the case when $V$ is a\qss
\emph{complex vector space of\dss dimension}\sss $1$\nnsp.\oss
It\sss turns out\sss that\sss in\sss this case $U\fff(V)$ has a canonical\sss structure
of a hyperbolic plane,\oss similar\sss to\sss the Klein model of the hyperbolic geometry,\oss
and\sss for $V\dff =\dff \cc$ the space $U\fff(V)\dff =\dff U\fff(\cc)$ 
can be identified with\sss the unit\sss disc\sss
$\uu
\dff =\dff 
\{\qff u\qff +\qff i\dff v\dff \mid\dff u^{\fff 2}\qff  +\qff  v^{\fff 2}\qff <\qff 1 \qff\}$\sss
in\sss the complex plane,\oss
which\sss is\sss nothing else but\sss the set\sss of\dss points of\dss the Klein model.\oss
It\sss is\dss also\sss the set\sss of\dss points of\dss the Poincar\'e model of\dss
hyperbolic geometry,\oss and\sss passing\sss from\sss the Klein model\sss
to\sss the Poincar\'e model\sss is\sss a key step in our approach\sss
to\sss the complex dilatation.\oss

\myuppar{The space $Q(\cc)$ of\dss quadratic forms on $\cc$\dnsp.}
For most\sss of\dss the paper we will\sss deal\sss with\sss the case $V\dff =\dff \cc$\nnsp.\oss
Our first\sss goal\dss is\dss to identify\sss $U\fff(\cc)$ with $\uu$\nnsp.\oss

Let\sss us denote by\dss  
$X\fff,\pff Y\dff\colon\dff \cc\toto\rr$\dss
the maps $X\dff(z)\off =\off \re z\dss and\dss
Y\dff(z)\off =\off \im z$\nnsp.\oss
As is well known,\qss $q$\dss is a quadratic form on\dss $\cc$\qss 
if and only if there exist\qss $a\fff,\pff b\fff,\pff c\dff\in\dff\rr$\qss
such that\qss 
$q(x\qff +\qff i\dff y)
\qff =\qff 
a\fff x^{\fff 2}\qff +\qff 2\fff b\fff x\fff y\qff +\qff c\fff y^{\fff 2}$\qss
for evey\sss $x\fff,\pff y\dff\in\dff \rr$\qss
or\halfff,\pss what is the same,\qss if
\[
\quad
q\off =\off a\fff X^{\dff 2}\qff +\qff2\fff  b\fff X\fff Y\qff +\qff c Y^{\dff 2}.
\] 
Clearly,\qss  $(X^{\dff 2},\off 2\fff X\fff Y,\off Y^{\dff 2})$\qss 
is basis of\dss $Q(\cc)$\dnsp,\qss
and the coefficients\dss $a\fff,\pff b\fff,\pff c$\dss 
are the coordinates of $q$ with respect to this basis.\oss
It is well known and easy to check that\sss 
$q\qff =\qff aX^{\dff 2}\qff +\qff 2\fff b X\fff Y\qff +\qff cY^{\dff 2}$\sss 
is a definite quadratic form\sss if\sss and only if\trs its\qss \emph{determinant}\qss
\[
\quad
D(q)\off =\off D(a\fff,\pff  b\fff,\pff  c)\off =\off ac\qff -\qff b^{\fff 2}
\]
is positive,\oss i.e.\oss $D(q)\qff =\qff ac\qff -\qff b^{\fff 2}\qff >\qff 0$\nnsp.\oss 
In other terms,\qss $[\dff q\dff]\dff\in\dff U\fff(\cc)$\qss if and only if\qss
$D(q)\qff >\qff 0$\dnsp.\oss
Let\sss us stress\sss the crucial fact is that\sss the determinant\sss $D$\sss 
is\dss itself\sss a quadratic form on\dss $Q(\cc)$\dnsp.\qss
This phenomenon is specific for dimension $2$\dnsp.\qss 

Another convenient\sss basis of\sss $Q(\cc)$ consists of\dss quadratic forms
\begin{equation}
\label{basis}
\quad
\mathbold{n}\off =\off
X^{\dff 2}\qff +\qff Y^{\dff 2},\pquad
\mathbold{r}\off =\off
X^{\dff 2}\qff -\qff Y^{\dff 2},\pquad
\mathbold{i}\off =\off
2\fff X\fff Y.
\end{equation}
The corresponding coordinates\sss $(t\fff,\pff r\fff,\pff s)$\sss
are related\sss to $(a\fff,\pff b\fff,\pff c)$
by\sss the equations
\[
\quad
s\off =\off b\dff,\quad
t\off =\off (a\qff +\qff c)/2\dff,\quad
r\off =\off (a\qff -\qff c)/2\dff.
\]
In\sss these coordinates $D$ has\sss the diagonal\sss form
$D(t\fff,\pff r\fff,\pff s)
\off =\off 
t^{\fff 2}\qff  -\qff  r^{\fff 2}\qff  -\qff  s^{\fff 2}$\dnsp.\oss
The basis\sss
$(\mathbold{n}\fff,\pff \mathbold{r}\fff,\pff \mathbold{i})$\sss 
has\sss the advantage of being closely related 
with\sss the structure of\dss $\cc$\nnsp.\oss
In fact\halfff,\oss
\[
\quad
\mathbold{n}\dff(z)
\off =\off 
|\dff z\dff|^{\fff 2}
\off =\off
z\fff\overline{z}\dff,\pquad
\mathbold{r}\dff(z)
\off =\off
\re z^{\fff 2},\pquad
\mathbold{i}\dff(z)
\off =\off
\im z^{\fff 2}
\]
for all\qss $z\dff\in\dff \cc$\nnsp.\oss
In particular\halfff,\qss 
$[\dff \mathbold{n}\dff]$\dss 
is the standard conformal structure on $\cc$\dnsp.\oss

\myuppar{A direct sum decomposition of $Q(\cc)$\dnsp.}
Now we are almost\sss ready\sss to identify\sss $U\fff(\cc)$ with $\uu$\nnsp.\qss
Let $\ccc$ be a copy of $\cc$\nnsp,\qss
which we will treat as a complex vector space.\qss 
Let us identify $Q(\cc)$ with $\rr \oplus \ccc$ by the map\dss
$(t\fff,\qff r\fff,\pff s)
\qff\longmapsto\qff (t\fff,\pff r\qff +\qff i\dff s)$\nnsp.\oss
Let\sss us set\sss
$z\qff =\qff r\qff +\qff i\dff s$
and consider\sss $(t\fff,\qff z)$
as\dss ``mixed''\sss coordinates on\sss
$Q(\cc)\qff =\qff \rr \oplus \ccc$\nnsp.\qss
In\sss the coordinates $(t\fff,\qff z)$\sss
the determinant $D$ takes the form\dss 
$D(t\fff,\pff z)
\qff =\qff
t^{\dff 2}\qff -\qff |\dff z\dff|^{\fff 2}$\dnsp.\qss
Hence\sss the space
$U\fff(\cc)$ of conformal structures on $\cc$
is defined by the homogeneous inequality\qss 
$t^{\dff 2}\qff -\qff |\dff z\dff|^{\fff 2}\qff >\qff 0$\nnsp,\oss
or\halfff,\oss equivalently,\oss by\sss the homogeneous inequality\sss
$|\dff z/t\dff|^{\fff 2}\qff <\qff 1$\nnsp.\oss

The plane\trs $\ccc\qff =\qff 0 \oplus \ccc$\trs 
in\trs $Q(\cc)\qff =\qff \rr \oplus \ccc$\qss
corresponds to a projective line in the projective plane\qss $\mathbb{P}\hnsp Q\fff(\cc)$\dnsp.\qss
The complement\sss $\mathbb{A}\fff(\cc)$\sss to this projective line
is an affine plane.\qss
Clearly,\dss $U\fff(\cc)\qff \subset\pff \mathbb{A}\fff(\cc)$\nnsp.\qss
We will\sss use\sss the map\sss
$(t\fff,\pff z)
\qff \longmapsto\qff
z/\halfff t
\qff\in\qff \ccc$
to identify\dss $\mathbb{A}\fff(\cc)$\dss with $\ccc$\nnsp.\qss
This identification\sss takes $[\dff\mathbold{n}\dff]$\sss to $0$\sss
and $U\fff(\cc)$\sss to\sss the unit\sss disc $\uu$ in $\ccc$\nnsp.\qss
\emph{This\sss is\sss the promised\sss identification of\dss
$U\fff(\cc)$ with $\uu$\dnsp.\qss}

Using\sss the copy\sss $\ccc$\sss of\sss $\cc$\sss allows us\sss 
to keep\sss the distinction between\sss the quadratic forms on\sss $\cc$\sss
represented\sss by elements of\sss $\ccc$
and complex numbers,\oss which are elements of\sss $\cc$\nnsp.\oss
By\sss this reason\sss the space $U\fff(\cc)$\sss is\dss identified with\sss
the unit\sss disc $\uu$ in\sss $\ccc$\nnsp,\oss not\sss in\sss $\cc$\nnsp.\oss

\myuppar{Diagonalization.}
For $a\fff,\pff c\dff\in\dff \rr$\sss let\sss
$q_{\fff a\fff,\dff c}\dff\in\dff Q(\cc)$\sss
be\sss the quadratic form\sss
$a\fff X^{\dff 2}\qff +\qff c Y^{\dff 2}$\dnsp.\oss
Obviously,\qss $q_{\fff a\fff,\dff c}$ is a definite form
if and only if $a\fff,\pff c$ are of the same sign.\oss
As\sss is\sss well\sss know,\qss every quadratic form on $\cc\dff =\dff \rr^{\dff 2}$
can be\sss turned\sss into\sss the form $q_{\fff a\fff,\dff c}$ 
for some $a\fff,\dff c$ by a rotation.\oss
Let\sss us state\sss this\qss \emph{diagonalization\sss theorem}\qss more precisely.\oss

To begin with,\qss
let $V\fff,\pff W$ be real\sss 
vector spaces and $L\colon V\toto W$ be a real linear map.\qss
For $p\dff \in\dff Q(W)$\sss let\sss the\dss \emph{pull-back}\dss
$L^*\dff p\dff\in\dff Q(V)$\sss be\sss the composition\sss 
$L^*\dff p\qff =\qff p\circ L$\nnsp.\oss
If\dss $U$ is another vector space and\qss $K\colon U\toto V$\qss is a linear map,\oss then\qss
$(L\circ\fff K)^*\qff =\qff K^*\circ\qff L^*$\dnsp.\oss

Let $V$ be a complex vector space.\qss
For $\tau\in \cc$\sss 
the\qss \emph{multiplication map}\dss
$m_{\dff \tau}\colon V\toto V$\qss 
is\sss defined by\sss
$m_{\dff \tau}(v)\qff =\qff \tau\dff v$\nnsp.\qss
Taking pull-backs\dss by\dss $m_{\dff \tau}$\dss leads to the induced map\sss
$m_{\dff \tau}^*\dff\colon\dff
Q(V)\ttoo Q(V)$\nnsp.\qss
If\dss $\tau\qff \neq\qff 0$\dnsp,\qss then $m_{\dff \tau}^*$ induces a map
$\mathbb{P}\halfff m_{\dff \tau}^{*} \dff\colon\dff 
\mathbb{P}\hnsp Q(V)\ttoo \mathbb{P}\hnsp Q(V)$\dnsp.\oss

Finally,\qss let\sss
$r_{\theta}\dff\colon\dff\rtwo\toto\rtwo$\qss 
be the counter-clockwise rotation 
by an angle $\theta\in\rr$ around the origin $(0\fff,\pff 0)\in\rtwo$\dnsp.\oss
As is well known,\qss the standard identification $\rtwo\toto\cc$ turns $r_{\theta}$
into the multiplication map
$\dis
m_{\dff \tau}\dff\colon\dff\cc\toto\cc$\nnsp,\oss
where\qss 
$\dis
\tau\qff =\qff e^{\fff i\dff\theta}
\qff =\qff 
\cos\theta\qff +\qff i\dff \sin\theta$\dnsp.\oss
Now we are ready\sss to state our\sss form of\sss the diagonalization\sss theorem.\oss

\mypar{Theorem.}{diag}  \emph{For every quadratic form\dss $q$\dss on\qss $\cc$\qss
there exist\qss $a\fff,\pff c\fff,\pff \theta\in\rrr$\qss 
such that}\vspace{3pt}
\begin{equation*}
\quad
q\off =\off m_{\dff \tau}^*\qff q_{\fff a\fff,\dff c}\qff,
\end{equation*}

\vspace{-12pt}\vspace{3pt}
\emph{where\qss $\tau\off =\off e^{\fff i\dff\theta}$\nnsp.\qff\oss
Moreover\halfff,\oss 
one can assume that\qss $|\dff a\dff|\qff \geq\qff  |\dff c\dff|$\dnsp.\oss}

\proof\qss Let\sss $q\fff\in\fff Q(\cc)\qff =\qff Q(\rr^{\dff 2})$\dnsp.\oss 
By the spectral theorem for symmetric operators in $\rtwo$ 
there is a basis\dss $(v\fff,\pff w)$\dss of\dss $\rtwo$\dss 
orthonormal with respect to\sss
$\mathbold{n}\dff(x\fff,\pff y)\qff =\qff x^2\qff +\qff y^2$\qss 
and orthogonal with respect to $q$\nnsp.\oss
After interchanging $v$ and $w$\nnsp,\oss if\dss necessary,\oss 
we may assume that\sss $(v\fff,\pff w)$ has
the same orientation as the standard basis.\pss
Then there is a rotation $r_{\theta}$ taking\sss 
$(v\fff,\pff w)$\sss to the standard basis of\dss $\rtwo$\nnsp.\dff\oss
If\dss $r_{\theta}$\dss is\sss such\sss a\sss rotation\halfff,\oss then\qss
$\dis
q\off =\off r_{\theta}^*\qff q_{\fff a\fff,\dff c}$\qss
for\dss some\qss $a\fff,\pff c\dff\in\dff \rr$\dnsp.
Since\oss
$r_{\pi/2}^*\qff q_{\fff a\fff,\dff c}
\qff =\qff  
q_{\fff c\fff,\dff a}$\oss
and\oss 
$r_{\theta}^*\qff\circ\qff r_{\pi/2}^*
\qff =\qff 
r_{\theta\qff +\qff \pi/2}^*$\nsp,\oss
one can assume that\qss $|\dff a\dff|\qff \geq\qff  |\dff c\dff|$\nnsp.\oss  \eproof

\myuppar{The conformal\sss class of\dss a diagonalized\sss form.}
Let\sss $a\fff,\pff c > 0$ and\sss $\tau\dff \in\dff \ccc$\nnsp,\dss $\tau\qff \neq\qff 0$\nnsp,\qss
and\sss let\sss us consider\sss the quadratic form
$q
\qff =\qff
m_{\dff \tau}^*\qff q_{\fff a\fff,\dff c}$\nsp.\oss
Our next\sss goal\sss is\sss to identify\sss the conformal\sss class
$[\fff q\dff]$ of\sss $q$ as an element\sss of\sss $\uu$\dnsp.\oss
Let\sss us consider first\sss the case
$q\qff =\qff q_{\fff a\fff,\dff c}$\nsp.\oss
The\sss 
$(t\fff,\pff r\fff,\pff s)$\dnsp-co\-or\-di\-nates\dss of\sss 
$q_{\fff a\fff,\dff c}$\sss
are\sss 
$(\dff (a\qff +\qff c)/2\fff,\pff (a\qff -\qff c)/2\fff,\pff 0 \dff)$\pss
and\dss hence\vspace{3pt}
\[
\quad
q_{\fff a\fff,\dff c}
\off =\off
\bigl(\dff (a\qff +\qff c)/2\fff,\pff (a\qff -\qff c)/2 \dff\bigr)
\qff\in \qff
\rr \oplus \ccc
\] 

\vspace{-12pt}\vspace{3pt}
after the identification\sss 
$Q(\cc)\qff =\qff \rr \oplus \ccc$\nnsp.\qff\oss
It\dss follows\dss that\vspace{3pt}
\[
\quad
[\fff q_{\fff a\fff,\dff c}\dff]
\off =\off 
\frac{\dff a\qff -\qff c\dff}{a\qff +\qff c}\qff \in\qff \ccc
\]

\vspace{-12pt}\vspace{3pt}
after the identification\sss 
$\mathbb{A}\fff(\cc)\qff =\qff \ccc$\nnsp.\qss
In order\sss to deal\sss with\sss the general\sss case\sss
$q
\qff =\qff
m_{\dff \tau}^*\qff q_{\fff a\fff,\dff c}$
we need\sss to study\sss the action of\dss the pull-backs
$m_{\dff \tau}^*$ on\sss the space of\dss quadratic forms $Q(\cc)$
and of\dss the induced\sss maps $\mathbb{P}\halfff m_{\dff \tau}^{*}$
on projective plane $\mathbb{P}\hnsp Q(\cc)$\nnsp.\oss

\mypar{Lemma.}{mult-on-forms} \emph{The map\oss
$\dis
m_{\dff \tau}^*\dff\colon\dff
Q(\cc)\ttoo Q(\cc)$\oss
respects the direct sum decomposition\qss
$Q(\cc)\qff =\qff \rr \oplus \ccc$\nnsp.\oss
It acts on the summand\dss $\rr$\dss as the multiplication by\qss
$|\dff \tau\dff|^{\fff 2}$\qss and on the summand\dss $\ccc$\dss
as the multiplication map\dss $m_{\dff \rho}$\nnsp,\oss
where\qss $\rho\qff =\qff \overline{\tau}^{\qff 2}$\nnsp.\oss}

\proof 
The proof\dss is\sss an almost\sss straightforward calculation.\qss
Obviously,\qss if\qss $z\dff\in\dff \cc$\nnsp,\qss
then\vspace{3pt}
\[
\quad
m_{\dff \tau}^*\dff \mathbold{n}\dff(z)
\off =\off 
|\dff \tau\dff|^{\fff 2}\trf |\dff z\dff|^{\fff 2}\dff,
\hspace*{1.2em}
m_{\dff \tau}^*\dff \mathbold{r}\dff(z)
\hspace*{0.15em}
\off =\off
\re \tau^{\dff 2}\dff z^{\fff 2}\dff,
\hspace*{1.2em}
\mbox{and}\hspace*{1.2em}
m_{\dff \tau}^*\dff \mathbold{i}\dff(z)
\hspace*{0.3em}
\off =\off
\im \tau^{\dff 2}\dff z^{\fff 2}\dff.
\]

\vspace{-12pt}\vspace{3pt}
If we present\qss $\tau^{\dff 2}$\qss in the form\qss 
$\tau^{\dff 2}\qff =\qff \alpha\qff +\qff i\dff \beta$\qss 
with\qss 
$\alpha\fff,\pff \beta\dff\in\dff \rr$\nnsp,\oss
then\vspace{3pt}
\[
\quad
\re \tau^{\dff 2}\dff z^{\fff 2}
\off =\off
\alpha\qff \re z^{\fff 2}\qff -\qff \beta\qff \im z^{\fff 2}
\hspace*{1.2em}
\mbox{and}\hspace*{1.2em}
\im \tau^{\dff 2}\dff z^{\fff 2}
\off =\off
\beta\qff \re z^{\fff 2}\qff +\qff \alpha\qff \im z^{\fff 2}
\]

\vspace{-12pt}\vspace{3pt}
for all\qss $z\dff\in\dff \cc$\nnsp.\oss 
It follows that\oss
$\dis
m_{\dff \tau}^*\dff \mathbold{n}
\off =\off 
|\dff \tau\dff|^{\fff 2}\qff \mathbold{n}$\oss and\vspace{3pt}
\begin{equation}
\label{mult-ri}
\quad
m_{\dff \tau}^*\dff \mathbold{r}
\off =\off
\alpha\dff \mathbold{r}\qff -\qff \beta\dff \mathbold{i}\dff,
\pquad\trf
m_{\dff \tau}^*\dff \mathbold{i}
\off =\off
\beta\fff \mathbold{r}\qff +\qff \alpha\trf \mathbold{i}\dff,
\end{equation}

\vspace{-12pt}\vspace{3pt}
and hence\dss 
$m_{\dff \tau}^*$\dss
leaves both summands $\rr$\nnsp,\qss $\ccc$ invariant\halfff,\oss
and acts on $\rr$ as the multiplication by\dss $|\dff \tau\dff|^{\fff 2}$\nnsp.\oss
Under the identification\qss
$Q(\cc)\qff =\qss \rr \oplus \ccc$\qss the forms\dss
$\mathbold{r}\fff,\pff \mathbold{i}$\dss correspond to\dss 
$1\fff,\pff i\dff\in\dff \ccc$\dss respectively,\oss
and the formulas\qss (\ref{mult-ri})\qss take the form\oss\vspace{3pt}
\begin{equation}
\label{mult-1i-ccc}
\quad
m_{\dff \tau}^*\dff(1)
\off =\off
\alpha\qff -\qff \beta\dff i\dff,
\pquad
m_{\dff \tau}^*\dff(i)
\off =\off
\beta\qff +\qff \alpha\dff i\dff.
\end{equation}

\vspace{-12pt}\vspace{3pt}
Therefore,\qss $m_{\dff \tau}^*$\dss acts on $\ccc$ 
as the multiplication by\qss 
$\alpha\qff -\qff i\dff\beta\qff =\qff \overline{\tau}^{\qff 2}$\nnsp.\oss  \eproof

\mypar{Theorem.}{action-on-mu} \emph{Suppose that\sss $\tau\dff\in\dff \cc$\sss
and\sss $\tau\qff \neq\qff 0$\dnsp.\oss 
Then the map\dss
$\mathbb{P}\halfff m_{\dff \tau}^{*}$\dss leaves both\dss
$\mathbb{A}\fff(\cc)$\dss and\dss $U\fff(\cc)$\dss invariant\halfff.\oss 
In particular\halfff,\qss
$\mathbb{P}\halfff m_{\dff \tau}^{*}$\qss induces a map}\dss\vspace{3pt}
\[
\quad
m_{\dff \tau}^{**}\qff\colon\qff
\mathbb{A}\fff(\cc)
\dff\ttoo\qff
\mathbb{A}\fff(\cc)\dff.
\]

\vspace{-12pt}\vspace{3pt}
\emph{The identification\dss
$\mathbb{A}\fff(\cc) 
\qff =\qff \ccc$\trs 
turns\dss $m_{\dff \tau}^{**}$\trs into\dss 
$m_{\dff \sigma}\dff\colon\dff \ccc\ttoo \ccc$\nnsp,\qss
where\dss
$\sigma\off =\off \overline{\tau}\fff/\tau$\nnsp.\oss}

\proof
Lemma\qss \ref{mult-on-forms}\qss implies that\dss 
$\mathbb{P}\halfff m_{\dff \tau}^{*}$\qss leaves the projective line
corresponding to the plane\dss $\ccc$\dss 
in\trs $Q(\cc)\qff =\qff \rr \oplus \ccc$\qss 
invariant and hence leaves\dss
$\mathbb{A}\fff(\cc)$\dss invariant\halfff.\oss
Lemma\qss \ref{mult-on-forms}\qss also implies
that if\oss
$m_{\dff \tau}^{*}\dff(t\fff,\pff z)\off =\off (t'\fff,\pff z')$\dnsp,\oss
then\pss
$\dis
|\dff t'\dff|\off =\off |\dff \tau\dff|^{\dff 2}\trf|\dff t'\dff|$\pss
and\pss
$\dis
|\dff z'\dff|\off =\off |\dff \tau\dff|^{\dff 2}\dff|\trf z'\dff|$\nnsp.\oss
It follows that\dss $m_{\dff \tau}^{*}$\dss leaves invariant\sss the cone in\dss $Q(\cc)$\dss 
defined by the homogeneous inequality\qss
$\dis
t^{\dff 2}\qff -\qff |\dff z\dff|^{\fff 2}\qff >\qff 0$\nnsp.\oss
Hence the induced map\qss 
$\mathbb{P}\halfff m_{\dff \tau}^{*}$\qss
leaves\dss $U\fff(\cc)$\dss invariant\halfff.\oss 

Let\sss 
$(t\fff,\pff z)\dff\in\dff \rr \oplus \ccc$\sss
be a representative of\dss 
$\mu\dff\in\dff \mathbb{A}\fff(\cc)$\nnsp.\oss 
Then $\mu\qff =\qff z/t$ under\sss the identification\sss
$\mathbb{A}\fff(\cc)\qff =\qff \ccc$\nnsp.\oss
Lemma\qss \ref{mult-on-forms}\qss implies\sss that\sss\vspace{3pt}
\[
\quad
m_{\dff \tau}^{*}\dff (t\fff,\pff z)
\off =\off
(\dff|\dff \tau\dff|^{\dff 2}\qff t\fff,\off \overline{\tau}^{\qff 2}\dff z\dff)
\]

\vspace{-12pt}\vspace{3pt}
and hence\sss that\sss 
$(\dff|\dff \tau\dff|^{\dff 2}\qff t\fff,\off \overline{\tau}^{\qff 2}\dff z\dff)$\qss
is a representative of\qss $m_{\dff \tau}^{**}\dff(\mu)$\nnsp.\oss
Hence\vspace{3pt} 
\[
\quad
m_{\dff \tau}^{**}\dff(\mu)
\off =\off
{\dis \overline{\tau}^{\qff 2}\dff z}\Bigl/
{\dis \dff|\dff \tau\dff|^{\dff 2}\qff t\dff}
\off =\off
{\dis \overline{\tau}^{\qff 2}}\bigl/
{\dis \dff|\dff \tau\dff|^{\dff 2}\dff}\qff
\bigl(\dff z/\fff t \dff\bigr)
\]

\vspace{-12pt}\vspace{3pt}
under\sss the identification\sss
$\mathbb{A}\fff(\cc)\qff =\qff \ccc$\nnsp.\oss
It follows that this identification 
turns\dss $m_{\dff \tau}^{**}$\dss into the multiplication map\dss 
$m_{\dff \sigma}\dff\colon\dff \ccc\ttoo \ccc$\nnsp,\oss where 
$\sigma
\qff =\qff
{\dis \overline{\tau}^{\qff 2}}/
{\dis \dff|\dff \tau\dff|^{\dff 2}\dff}
\off =\off
{\dis \overline{\tau}^{\qff 2}}/
{\dis \dff\tau\dff \overline{\tau}\dff}
\off =\off
\overline{\tau}/\fff\tau$\nnsp.\oss  \eproof

\mypar{Corollary.}{action-on-q} \emph{If\qss $q\dff\in\dff Q(\cc)$\qss 
is\dss a\dss definite\dss form\dss 
and\oss 
$\tau\dff\in\dff \cc$\nnsp,\oss $\tau\qff \neq\qff 0$\dnsp,\qff\oss 
then}\oss\vspace{1.5pt} 
\[
\quad
[\dff m_{\dff \tau}^*\qff q\dff]
\off =\off
(\fff\overline{\tau}\fff/\tau\fff)\qff[\dff q\dff]\dff\in\dff\ccc\dff.
\hspace*{1.5em}\mbox{\eproof}
\]

\vspace{-12pt}\vspace{1.5pt}
\mypar{Corollary.}{uniqueness} \emph{There is 
a unique conformal structure on\dss $\cc$\dss
fixed by all maps\qss 
$\mathbb{P}\halfff m_{\dff \tau}^{*}$\qss 
with\qss
$\tau\qff \neq\qff 0$\dnsp.\qff\oss
It is equal to the conformal class\pss $[\dff\mathbold{n}\dff]$\pss of\pss $\mathbold{n}$\nnsp.\oss}

\proof Theorem\qss \ref{action-on-mu}\qss implies 
that a conformal structure on\dss $\cc$\dss is
fixed by all maps\dss 
$\mathbb{P}\halfff m_{\dff \tau}^{*}$\dss 
with\dss
$\tau\qff \neq\qff 0$\dss 
if and only if  
as a point of\dss $\mathbb{U}$\dss it is
fixed by the multiplication by all\dss $\sigma\dff\in\dff\ccc$\dss of the form\qss 
$\sigma\qff =\qff \overline{\tau}\fff/\tau$\dss with\qss $\tau\qff \neq\qff 0$\dnsp.\oss
But\dss $\sigma$\dss has such form if and only if\qss
$|\dff\sigma\dff|\qff =\qff 1$\nnsp,\oss
and there is exactly one point of\dss $\mathbb{U}$\dss fixed by the multiplication
by all such\dss $\sigma$\nnsp,\oss 
namely,\oss the point $0$\dnsp.\oss
Since $0$ corresponds to\dss $[\dff\mathbold{n}\dff]$\nnsp,\oss
this proves the corollary.\oss  \eproof

\mypar{Theorem.}{conformal-class-diag} \emph{Let\pss 
$a\fff,\pff c\qff\in\qff \rr$\pss and\pss $\tau\qff\in\qff \cc$\nnsp.\oss
Suppose that\qss
$a\fff,\pff c$\pss are\dss positive\dss
and\pss $\tau\qff \neq\qff 0$\dnsp.\oss
Let\sss
$q
\off =\off
m_{\dff \tau}^*\qff q_{\fff a\fff,\dff c}$\nsp.\oss
Then under\sss the identification\sss
$\mathbb{A}\fff(\cc)\qff =\qff \ccc$\sss
the conformal\sss class}
\[
\quad
[\dff q\dff]
\off =\off
\bigl(\dff\overline{\tau}\fff/\tau\dff\bigr)\qff
\frac{\dff a\qff -\qff c\dff}{\dff a\qff +\qff c\dff} 
\pff \in\pff \ccc
\qff.
\]

\vspace{-12pt}\vspace{-6pt}
\proof
As we already saw,\qss
$\dis
[\fff q_{\fff a\fff,\dff c}\dff]
\off =\off 
\frac{\dff a\qff -\qff c\dff}{a\qff +\qff c}\qff \in\qff \ccc$\nnsp.\oss
By\sss the definition of\sss $m_{\dff \tau}^{**}$\dss
we have
\[
\quad
[\fff m_{\dff \tau}^{*}\qff q_{\fff a\fff,\dff c}\dff]
\off =\off
m_{\dff \tau}^{**}\qff [\fff q_{\fff a\fff,\dff c}\dff]
\off =\off
m_{\dff \tau}^{**}\left(\frac{\dff a\qff -\qff c\dff}{a\qff +\qff c}\right).
\]
It remains to apply Theorem\qss \ref{action-on-mu}.\oss  \eproof

\mysection{The\qss Poincar\'e\qss invariant\qss of\qss a\qss quadratic\qss form}{p-invariant}

\vspace{6pt}
\myuppar{From the Klein model to the Poincar\'e model.} 
The open unit disc\sss $\uu\dff\subset\dff \ccc$\sss 
is the set of points of the Klein model of hyperbolic geometry.\qss
It is also the set of points of the Poincar\'e model of hyperbolic geometry.\qss
There is a canonical isomorphism between these two structures 
on\sss $\uu$\dnsp.\oss
In order to describe it\halfff,\qss
let us use the upper hemisphere\vspace{2pt} 
\[
\quad
S^{\fff 2}_{+}\off =\off 
\{\qff(u\qff +\qff i\dff v\fff,\pff w) \qff\mid\qff 
u^{\fff 2}\qff +\qff v^{\fff 2}\qff +\qff w^{\fff 2}\qff =\qff 1\fff,\qff\off w\qff >\qff 0 \qff\} 
\off\subset\off \ccc\times\rr
\] 

\vspace{-12pt}\vspace{2pt}
as an intermediary between the Klein and the Poincar\'e models.\oss 

Let $\pi$ be the orthogonal projection of the upper hemisphere 
$S^{\fff 2}_{+}$ to the equatorial disc 
$\uu\off =\off
\{\qff (u\qff +\qff i\dff v)\dff\in\dff\ccc \qff\mid\qff u^2\qff +\qff v^2\qff<\qff 1 \qff\}$\nsp.\oss
In other terms,\dss 
$\pi\trf(u\qff +\qff i\dff v\fff,\pff w)
\off =\off 
u\qff +\qff i\dff v
\qff\in\qff \ccc$\nnsp.\qss
Let $S$ be the stereographic projection of the upper hemisphere $S^{\fff 2}_{+}$ 
to the equatorial disc from the point\qss 
$s\qff =\qff (0\qff +\qff i\dff 0\fff,\pff -\qff 1)$\dnsp.\oss 
As\sss is\sss well\sss known,\qss the map 
$\Omega\qff =\qff S\circ\pi^{\qff-\dff 1}$ 
transforms the Klein structure of the hyperbolic plane on $\uu$ 
into the Poincar\'e one.\oss

\myuppar{The Poincar\'e invariant of quadratic forms on $\cc$\dnsp.}
For every definite quadratic form\dss $q$\dss on\dss $\cc$\dss we define
its\qss \emph{Poincar\'e invariant}\dss $P\fff(q)$ as\sss the image\sss 
$P\fff(q)
\qff =\qff 
\Omega\dff\bigl(\fff[\fff q\dff]\fff\bigr)
\dff\in\dff \uu$
of its conformal class $[\fff q\dff]$ under the map $\Omega $\nnsp.\oss
Our next\sss goal\sss is\dss to compute\sss $P\fff(q)$ for 
$q
\off =\off
m_{\dff \tau}^*\qff q_{\fff a\fff,\dff c}$\nsp.\oss
To begin with,\qss observe\sss that
\begin{equation}
\label{p-rotations}
\quad
P\fff (\fff m_{\dff \tau}^*\qff q\fff )
\off =\off
(\overline{\tau}\fff/\tau)\qff P\fff (\fff q\fff )
\end{equation}
for every $q\dff\in\dff Q(\cc)$
and\sss $\tau\dff\in\dff \cc$\nnsp,\dss $\tau\dff \neq\dff 0$\nnsp.\oss
Indeed,\qss the map $\Omega$ 
obviously commutes with the multiplication maps 
$m_{\dff \sigma}$ for 
$|\dff \sigma\dff |\qff =\qff 1$\nnsp,\qss
and\sss hence\qss (\ref{p-rotations})\qss follows from\sss Corollary\qss \ref{action-on-q}.\qss
We need also an explicit\sss description of\dss the map $\Omega$\nnsp.\oss

\mypar{Lemma.}{k-to-p} 
\emph{Suppose\dss that\qss $\gamma\in\rr$\nnsp,\oss 
$\sigma\in\ccc$\oss
and\oss $\gamma > 0$\dnsp,\oss $|\dff\sigma\dff|\qff =\qff 1$\dnsp.\oss}\oss
\[
\quad
\hspace*{1.5em}\mbox{\emph{If}}\hspace*{1em}
p
\off =\off
\frac{1\qff -\qff \gamma^{\dff 2}}
{\dff 1\qff +\qff \gamma^{\dff 2}\dff}\qff\sigma\dff,
\hspace*{1.5em}\mbox{ \emph{then} }\hspace*{1.5em}
\Omega\dff(\fff p)
\off =\off
\frac{1\qff -\qff \gamma}
{\dff 1\qff +\qff \gamma\dff}\qff\sigma\qff.
\]

\vspace{-12pt}
\proof The counter-clockwise rotation of\dss $\ccc\times\rr$\dss 
along the vertical axis\dss $0\times\rr$\dss 
by the angle\qss $-\dff\arg \sigma$\qss reduces the proof 
to the case\qss $\sigma\qff =\qff 1$\nnsp.\oss
In this case everything happens in the vertical coordinate plane\qss $v\qff =\qff 0$\dnsp,\oss
and we can ignore the\dss $v$\dnsp-coordinate and write\vspace{3pt}
\[
\quad
p
\off =\off
\left(\frac{1\qff -\qff \gamma^{\dff 2}}
{\dff 1\qff +\qff \gamma^{\dff 2}\dff}\dff,\off 0\right).
\]

\vspace{-12pt}\vspace{3pt}
By using the well known identity\off\oss\vspace{3pt}
\[
\quad
\left(\frac{1\qff -\qff \gamma^{\dff 2}}
{\dff 1\qff +\qff \gamma^{\dff 2}\dff}\right)^2
\qff +\qff
\left(\frac{2\dff\gamma}
{\dff 1\qff +\qff \gamma^{\dff 2}\dff}\right)^2
\off =\off 1
\]

\vspace{-12pt}\vspace{3pt}
and the assumption\qss $\gamma\qff >\qff 0$\dnsp,\oss 
we see that\oss\vspace{2pt}
\[
\quad
\pi^{\dff -\dff 1}\fff(p)
\off =\off
\left(\frac{1\qff -\qff \gamma^{\dff 2}}
{\dff 1\qff +\qff \gamma^{\dff 2}\dff}\dff,\off 
\frac{2\dff\gamma}
{\dff 1\qff +\qff \gamma^{\dff 2}\dff}\right).
\]

\vspace{-12pt}\vspace{2pt}
We need to prove that the stereographic projection $S$ takes this point to\vspace{2pt}
\[
\quad
\left(\frac{1\qff -\qff \gamma}
{\dff 1\qff +\qff \gamma\dff}\dff,\off 0\right).
\]

\vspace{-12pt}\vspace{2pt}
By the definition of $S$ this means that the points\vspace{1.5pt}
\[
\quad
(0\fff,\pff  -\qff 1),
\hspace*{1.5em}
\left(\frac{1\qff -\qff \gamma}
{\dff 1\qff +\qff \gamma\dff}\dff,\off 0\right),
\hspace*{1.5em}\mbox{and}\hspace*{1.5em}
\left(\frac{1\qff -\qff \gamma^{\dff 2}}
{\qff 1\qff +\qff \gamma^{\dff 2}\dff}\dff,\off 
\frac{2\dff\gamma}
{\dff 1\qff +\qff \gamma^{\dff 2}\dff}\right)
\]

\vspace{-12pt}\vspace{1.5pt}
are contained in the same line.\oss
This happens if and only if the vectors\vspace{1.5pt}
\[
\quad
\left(\frac{1\qff -\qff \gamma}
{\dff 1\qff +\qff \gamma\dff}\dff,\off 0\right)
\off -\off
(0\fff,\pff  -\qff 1)
\off =\off
\left(\frac{1\qff -\qff \gamma}
{\dff 1\qff +\qff \gamma\dff}\dff,\off 1\right)
\quad\mbox{and}
\]

\vspace{-24pt}
\[
\quad
\left(\frac{1\qff -\qff \gamma^{\dff 2}}
{\qff 1\qff +\qff \gamma^{\dff 2}\dff}\dff,\off 
\frac{2\dff\gamma}
{\qff 1\qff +\qff \gamma^{\dff 2}\qff}\right)
\off -\off
(0\fff,\pff  -\qff 1)
\off =\off
\left(\frac{1\qff -\qff \gamma^{\dff 2}}
{\qff 1\qff +\qff \gamma^{\dff 2}\dff}\dff,\off 
\frac{(\dff 1\qff +\qff \gamma\dff)^2}
{\qff 1\qff +\qff \gamma^{\dff 2}\dff}\qff\right)
\]

\vspace{-12pt}\vspace{1.5pt}
are proportional\halfff,\oss i.e.\qss if and only if\vspace{1.5pt}
\[
\quad
\frac{1\qff -\qff \gamma^{\dff 2}}
{\qff 1\qff +\qff \gamma^{\dff 2}\dff}
\left/
\frac{1\qff -\qff \gamma}
{\dff 1\qff +\qff \gamma\dff}
\off =\off
\frac{(\dff 1\qff +\qff \gamma\dff)^2}
{\dff 1\qff +\qff \gamma^{\dff 2}\dff}
\right.\dff.
\]

\vspace{-12pt}\vspace{1.53pt}
A straightforward verification of this identity completes the proof\halfff.\oss  \eproof

\mypar{Theorem.}{form-to-p} \emph{If\dss the\dss form\oss  
$\dis
q\off =\off m_{\dff \tau}^*\qff q_{\fff a\fff,\dff c}$\qss 
is\dss definite,\qff\oss 
then\qss $c/a$\qss is\dss positive\dss and}\oss\vspace{1.5pt}
\[
\quad
P\fff(\fff q\fff )
\off =\off
\frac{1\qff -\qff \gamma}{\dff 1\qff  +\qff \gamma\dff}\qff (\overline{\tau}\fff/\tau),
\]

\vspace{-12pt}\vspace{1.5pt}
\emph{where\dss $\gamma$\dss is\dss the\dss positive\dss square\dss root\qss $\sqrt{c/a\dff}$\nsp.\oss}

\proof If\qss $m_{\dff \tau}^*\qff q_{\fff a\fff,\dff c}$\qss is definite,\oss
then\sss $q_{\fff a\fff,\dff c}$\sss is also definite,\oss 
and hence the real numbers\sss 
$a\fff,\pff c$\sss are of the same sign.\oss
Therefore\qss $c/a$\qss is positive.\oss
By\dss Theorem\qss \ref{conformal-class-diag}\qss\vspace{3pt}  
\[
\quad
[\fff q\dff]
\qff\off =\qff\off
\frac{\dff a\qff -\qff c\dff}{\dff a\qff +\qff c\dff}\off (\overline{\tau}\fff/\tau)
\off =\off
\frac{1\qff -\qff c/a}{\dff 1\qff +\qff c/a\dff}\off (\overline{\tau}\fff/\tau)
\qff\off =\qff\off
\frac{1\qff -\qff \gamma^{\dff 2}}
{\dff 1\qff +\qff \gamma^{\dff 2}\dff}\off (\overline{\tau}\fff/\tau)
\dff.
\] 

\vspace{-12pt}\vspace{3pt}
An application of\qss Lemma\qss \ref{k-to-p}\qss completes the proof\halfff.\oss  \eproof

\myuppar{Remark.}
Lemma\qss \ref{k-to-p}\qss leads\sss to a simple formula for\sss the map $\Omega^{-\dff 1}$\dnsp.\qss
Given\qss $\mu\dff\in\dff \uu$\dnsp,\qss 
let\sss
$\sigma
\qff =\qff 
\mu\fff/\fff |\dff \mu\dff|$\sss
and\sss
$\gamma
\qff =\qff
(1\qff -\qff |\dff \mu\dff|\fff)/(1\qff +\qff |\dff \mu\dff|\fff)$\nnsp.\qss
Then\sss $|\dff\sigma\dff|\qff =\qff 1$\nnsp,\dss $\gamma > 0$\nnsp,\qss\vspace{3pt}
\[
\quad
|\dff \mu\dff|
\off =\off
\frac{1\qff -\qff \gamma}
{\dff 1\qff +\qff \gamma\dff}
\dff,
\hspace*{1.2em}
\mu
\off =\off
\frac{1\qff -\qff \gamma}
{\dff 1\qff +\qff \gamma\dff}\qff\sigma
\dff,\hspace*{1.2em}
\frac{1\qff -\qff \gamma^{\dff 2}}
{\dff 1\qff +\qff \gamma^{\dff 2}\dff}
\dff\off =\dff\off
\frac{2\dff |\dff \mu\dff|}
{\dff 1\qff +\qff |\dff \mu\dff|^{\dff 2}\dff}
\qff,
\hspace*{1.2em}
\mbox{and\dss hence}
\]

\vspace{-24.5pt}
\[
\quad
\Omega^{-\dff 1}(\mu)
\off =\off
\frac{1\qff -\qff \gamma^{\dff 2}}
{\dff 1\qff +\qff \gamma^{\dff 2}\dff}\qff\sigma
\off =\off
\frac{2\dff |\dff \mu\dff|}
{\dff 1\qff +\qff |\dff \mu\dff|^{\dff 2}\dff}\qff\sigma
\off =\off
\frac{2\dff \mu}{\dff 1\qff +\qff |\dff \mu\dff|^{\dff 2}\dff}\qff.
\]

\mysection{The\qss dilatation\qss of\qss real\qss linear\qss maps\qss $\cc\toto\cc$}{linear maps}

\vspace{6pt}
\myuppar{The Poincar\'e dilatation of a real linear map $\cc\rightarrow\cc$\nnsp.} 
For an invertible real linear map\sss $T\dff \colon\dff \cc \toto\cc$\sss 
the pull-back\sss 
$T^{\dff *}\dff [\dff \mathbold{n}\dff]\qff\in\qff \uu$\sss
is a natural measure of distortion of the standard conformal structure 
$[\dff \mathbold{n}\dff ]$ by $T$\dnsp.\qss
Let\sss us define\qss \emph{Poincar\'e conformal dilatation}\dss $\pi_{\trf T}$ as\vspace{3pt} 
\[
\quad
\pi_{\trf T}
\off =\off
\Omega\qff\bigl(\dff T^{\dff *}\dff [\dff \mathbold{n}\dff] \dff\bigr)
\dff.
\] 

\vspace{-12pt}\vspace{3pt}
In order to relate\sss the Poincar\'e conformal dilatation $\pi_{\trf T}$ 
with the classical complex dilatation $\mu_{\dff T}$ of\dss $T$\nnsp,\oss
let us consider the image $E$ 
of the unit circle in $\cc$ under $T$\nnsp.\oss
As is well known,\qss 
$E$ is an ellipse with the center at the origin.

\mypar{Theorem.}{geom-interpretation} (i)\oss  
\emph{Let\dss $D$\dss be the ratio of the major axis of\qss $E$\dss 
to the minor axis.\oss Then}\vspace{1.5pt}
\[
\quad
|\dff\pi_{\trf T}\dff|\off =\off \frac{\dff D\qff -\qff 1\dff}{D\qff +\qff 1}\dff.
\]

\vspace{-12pt}\vspace{1.5pt}
(ii)\oss 
\emph{If\qss $\alpha$\qss is an angle between $\rr\dff\subset\dff\cc$ and
a direction of the maximal distortion of\qss $T$\dnsp,\oss then}\vspace{1.5pt}
\[
\quad
\arg\dff \pi_{\trf T}\off =\off 2\halfff\alpha\fff.
\]

\vspace{-12pt}\vspace{1.5pt} 
\proof Let\sss
$q\qff =\qff T^{\dff *}\dff\mathbold{n}$\nnsp.\qss
Then\sss $\pi_{\dff T}\qff =\qff P\fff(q)$\sss
and\sss
$q\fff(z)
\qff =\qff
\mathbold{n}\dff(\dff T\fff(z)\dff)
\qff =\qff
\|\dff T\fff(z)\dff \|^{\fff 2}$\sss
for every\qss $z\dff\in\dff\cc$\nnsp.\oss
By\dss Theorem\qss \ref{diag},\dss 
$q\qff =\qff m_{\dff \tau}^*\pff q_{\fff a\fff,\dff c}$\sss
for some\qss $a\fff,\pff c\dff\in\dff \rr$\qss and\qss $\tau\dff\in\dff \cc$\qss
such that\qss $|\dff\tau\dff|\qff =\qff 1$\qss and\qss
$|\dff a\dff|\qff \geq\qff  |\dff c\dff|$\dnsp.\oss
Since $q$ is positive definite together with $\mathbold{n}$\nnsp,\qss 
we may assume\sss that\sss $a\qff \geq\qff c\qff >\qff 0$\nnsp.\qss
Let\oss $\gamma\qff =\qff \sqrt{c/a}$\nsp.\qss
For every\sss $z\dff\in\dff\cc$\vspace{3pt}
\[
\quad
\|\dff T(z)\dff \|^{\fff 2}
\off =\off
q\fff(z)
\off =\off
m_{\dff \tau}^*\qff q_{\fff a\fff,\dff c}\dff(z)
\off =\off
q_{\fff a\fff,\dff c}\fff(\dff m_{\dff \tau}\dff(z)\dff)
\off =\off
q_{\fff a\fff,\dff c}\fff(\dff \tau\dff z\dff)
\dff.
\]

\vspace{-12pt}\vspace{3pt}
At\sss the same\sss time\sss
$\|\dff z\dff\|^{\fff 2}
\off =\off
\|\dff \tau\dff z\dff\|^{\fff 2}$\oss
because\oss $|\dff\tau\dff |\off =\off 1$\dnsp.\qss
It follows that the square of the distortion of\dss $T$\sss 
in the direction of a non-zero vector\sss 
$z\dff\in\dff \cc\qff =\qff \rtwo$\sss is\sss equal to\vspace{3pt}
\[
\quad
\|\dff T(z)\dff \|^2 
\bigl/\dff 
\|\dff z\dff \|^{\fff 2}
\off =\off
\|\dff q_{\fff a\fff,\dff c}\fff(\dff \tau\dff z\dff)\dff \|^{\fff 2}
\bigl/\dff
\|\dff z\dff \|^{\fff 2}
\off =\off
\|\dff q_{\fff a\fff,\dff c}\fff(\dff \tau\dff z\dff)\dff \|^{\fff 2}
\bigl/\dff
\|\dff \tau\dff z\dff \|^{\fff 2}
\qff.
\]

\vspace{-12pt}\vspace{3pt}
Since\sss $a\qff \geq\qff c$\dnsp,\qss 
the ratio\sss
$\|\dff q_{\fff a\fff,\dff c}\fff(\dff w\dff)\dff \|^2
\bigl/\dff
\|\dff w\dff \|^2$
achieves its maximum\dss $a$\dss when\sss $w\dff\in\dff \rr$\nnsp,\qss
and its minimum $c$ when\dss $w\dff\in\dff i\dff\rr$\nnsp.\qss
Hence\sss
$\|\dff T(z)\dff \|
\bigl/\dff
\|\dff z\dff \|$\sss
achieves its maximum\dss $\sqrt{a}$\dss when\qss
$\tau\dff z\dff\in\dff \rr$\nnsp,\oss
and its minimum\dss $\sqrt{c}$\dss when\qss 
$\tau\dff z\dff\in\dff i\dff\rr$\nnsp.\oss 

The length of the major axis of the ellipsis $E$ 
is equal to the the maximal distortion
of\qss $T$\dnsp,\oss i.e.\qss to\qss $\sqrt{a}$\nsp,\oss
and the length of its minor axis 
is equal to the the minimal distortion,\oss 
i.e.\qss to\qss $\sqrt{c}$\nsp.\oss
It follows that the ratio $D$ of the axes of\sss $E$ is equal to 
$\gamma^{\dff -\dff 1}\qff =\qff \sqrt{a/c}$\nsp.\oss
Since\qss $|\dff\overline{\tau}\fff/\tau\dff|\qff =\qff 1$\dnsp,\oss
this fact together with\dss Theorem\qss \ref{form-to-p}\qss 
implies the part\qss (i)\qss of the\sss theorem.\qss

The distortion of\dss $T$\dss in the direction of vector\dss $z$\dss 
is maximal when\sss 
$w\qff =\qff \tau\dff z\dff\in\dff \rr\dff\subset\dff \cc$\nnsp,\qss
or\halfff,\qss what\dss is\dss the same,\qss when\sss
$\arg (\tau\dff \dff z)\qff =\qff 0$\qss or\qss $\pi$\nnsp.\qss
Equivalently,\qss the distortion is maximal\sss when\sss
$\arg z\qff =\qff -\dff \arg\fff \tau$\qss or\qss 
$-\dff \arg\fff \tau\qff +\qff \pi$\nnsp.\qss
Therefore\sss $\alpha\qff =\qff -\dff \arg\fff \tau$\sss or\sss
$-\dff \arg\fff \tau\qff +\qff \pi$\sss and hence\sss
$2\fff\alpha\qff =\qff -\qff 2\fff\arg\fff \tau$\nnsp.\qss
On the other hand,\qss 
$2\fff\arg\fff \tau
\qff =\qff
-\qff \arg\dff (\overline{\tau}\fff/\tau)$\qss and\dss
Theorem\qss \ref{form-to-p}\qss implies that\qss
$\arg\dff (\overline{\tau}\fff/\tau)
\qff =\qff
\arg\dff \pi_{\dff T}$\nnsp.\oss
Hence $2\fff\alpha\qff =\qff \arg\dff \pi_{\dff T}$\nnsp.\qss
This proves the part\qss (ii).\oss  \eproof

\mypar{Main\dss Theorem.}{p-equal-classical} 
\emph{The Poincar\'e conformal dilatation of
a real linear orientation-preserving automorphism\qss $T\colon \cc\toto\cc$\qss
is equal to its classical complex dilatation.\oss}

\proof By Theorem\qss \ref{geom-interpretation},\qss the Poincar\'e conformal dilatation
satisfies properties\qss (i)\qss and\qss (ii)\qss from the Introduction\qss
(or\halfff,\qss rather\halfff,\qss their analogues for real linear maps)\qss
characterizing the classical complex dilatation.\oss  \eproof

\mysection{The action of real linear maps\qss $\cc\ttoo \cc$}{action}

\vspace{6pt}
\myuppar{Automorphisms induced by a real linear invertible maps\qss $\cc\toto \cc$\nnsp.}
Let\qss $u\dff\colon\dff \cc\toto \cc$\qss is a real linear invertible map.\oss
Then\sss the pull-back map\sss
$u^*\qff \colon\qff Q(\cc)\ttoo Q(\cc)$
is an automorphism of the real vector space\dss $Q(\cc)$\dss and hence induces an automorphism
\[
\quad
\mathbb{P}\fff u^*\qff \colon\qff
\mathbb{P}\hnsp Q(\cc)\ttoo \mathbb{P}\hnsp Q(\cc)
\]
of the projective plane $\mathbb{P}\hnsp Q(\cc)$\dnsp.\qss
Clearly\halfff,\oss a quadratic form $q\dff\in\dff Q(\cc)$
is definite if and only if the pull-back $u^*\fff q$\sss is definite,\oss
and hence $\mathbb{P}\fff u^*$\sss leaves the space 
$U\fff(\cc)$ 
of conformal structures on\dss $\cc$\dss invariant\halfff.\oss
After\sss the identification of $Q(\cc)$ with $\rr \oplus \ccc$
and $\mathbb{A}\fff(\cc)$ with $\ccc$ the space $U\fff(\cc)$
turns into the open unit disc $\uu$ in the complex plane $\ccc$\nnsp,\qss
and the automorphism\sss $\mathbb{P}\fff u^*$\sss 
turns into an automorphism of the projective plane\sss
$\mathbb{P}(\rr \oplus \ccc)$\sss
leaving\dss $\uu$\dss invariant\halfff.\oss

The structure of the Klein model of the hyperbolic plane on $\uu$
is defined in terms of
$\uu$ and the geometry of the projective plane\sss $\mathbb{P}(\rr \oplus \ccc)$\dnsp.\qss
Therefore this structure is invariant under any automorphism of this
projective plane leaving $\uu$ invariant\halfff.\oss
In particular\halfff,\qss it is invariant under\sss $\mathbb{P}\fff u^*$\nnsp.\qss
In other terms,\dss
$\mathbb{P}\fff u^*$ is an automorphism of the Klein model.\oss
Passing to the Poincar\'e unit\sss disc model\sss of the hyperbolic plane turns $\mathbb{P}\fff u^*$
into an automorphism of the Poincar\'e model.\qss
In more details,\oss using the map $\Omega$ to pass to the Poincar\'e model
turns\sss $\mathbb{P}\fff u^*$ into the map\qss
(denoted by $u^*$ again)
\[
\quad
u^*
\off =\off
\Omega\dff\circ\dff \mathbb{P}\fff u^*\circ\dff \Omega^{-\dff 1}
\qff \colon\qff
\uu \ttoo \uu
\]
and this map is an automorphism of the Poincar\'e unit disc model.\oss 
It\sss is\sss easy\sss to see\sss that\sss if\dss $u$\dss is orientation preserving\halfff,\oss
then\sss 
$u^*\dff \colon\dff \uu\toto \uu$\qss is also orientation-preserving.\qss

\mypar{Lemma.}{action-on-pi}
\emph{If\oss $T\dff \colon\dff \cc\ttoo \cc$\oss is a real linear invertible map,\qff\oss
then\sss
$\pi_{\trf T\dff\circ\dff u}
\off =\off
u^*\dff(\pi_{\trf T})$\nnsp.\qss}

\proof
Obviously\halfff,\dss
$\dis
(\fff T\circ u\fff)^*\dff \mathbold{n}
\off =\off
(\fff u^*\circ T^{\dff *}\fff)\dff (\mathbold{n})
\off =\off
u^*\dff (\fff T^{\dff *}\dff \mathbold{n}\fff)$\dss
and hence\vspace{3pt}
\[
\quad
[\dff (\fff T\circ u\fff)^*\dff \mathbold{n}\dff]
\off =\off
\mathbb{P}\fff u^*\dff \bigl(\dff [\dff T^{\dff *}\dff \mathbold{n}\dff] \dff\bigr)\dff.
\]

\vspace{-12pt}\vspace{3pt}
By applying\sss $\Omega$\sss to the both part of this equality
we see that\vspace{3pt}
\[
\quad
\pi_{\trf T\dff\circ\dff u}
\off =\off
\Omega\dff \bigl(\dff [\dff (T\circ u)^*\dff \mathbold{n}\dff] \dff\bigr)
\off =\off
\Omega\dff\circ\dff \mathbb{P}\fff u^*\circ\dff \Omega^{-\dff 1}
\bigl(\dff \Omega\dff [\dff T^*\dff \mathbold{n}\dff] \dff\bigr)
\off =\off
u^*\dff(\pi_{\trf T})\dff.
\hspace*{1.5em}\mbox{ \eproof }
\]

\vspace{-12pt}\vspace{3pt}
\mypar{Theorem.}{action-on-beltrami}
\emph{If\dss $u$\sss is orientation-preserving\halfff,\qff\oss
then\sss $u^*$\sss is complex analytic and}\vspace{0pt}
\[
\quad
\mu_{\dff T\dff\circ\dff u}
\off =\off
u^*\dff(\mu_{\dff T})
\]

\vspace{-12pt}\vspace{0pt}
\emph{for every real linear orientation-preserving invertible map\oss 
$T\dff \colon\dff \cc\ttoo \cc$\nnsp.\oss}

\proof
If\sss $u$ is orientation-preserving\halfff,\qss
then $u^*$ is\sss also orientation-preserving.\qss
Since every orientation-preserving automorphism of the Poincar\'e model
is complex analytic,\oss it follows that\sss $u^*$ is complex analytic.\oss
Moreover\halfff,\qss
if\sss $u$ is orientation-preserving\halfff,\qss
then $T\dff\circ\dff u$ is orientation-preserving together with $T$\dnsp.\qss
Therefore,\qss Theorem\qss \ref{p-equal-classical}\qss implies that\sss
$\mu_{\dff T\dff\circ\dff u}
\qff =\qff
\pi_{\dff T\dff\circ\dff u}$\sss and\sss
$\mu_{\dff T}
\qff =\qff
\pi_{\dff T}$\nsp.\qss
It remains to apply Lemma\qss \ref{action-on-pi}.\oss  \eproof

\mysection{Complex vector spaces of dimension 1}{one-dim-over-cc}

\vspace{6pt}
\myuppar{Quadratic forms associated with non-zero vectors.} 
The goal\sss of\sss this section is\sss to extend\sss the notion of\dss the\sss 
Poincar\'e conformal dilatation\sss from $\cc$\sss 
to complex vector spaces $V$ of\sss dimension $1$\nnsp.\qss
Let\sss us start\sss with elementary observations about\sss such vector spaces.\qss
A basis of\sss $V$ over $\cc$
is just a non-zero vector\qss $v\in V $\dnsp.\qss
Hence,\qss for any non-zero\qss $v\in V$\qss there is a unique 
isomorphism of complex vector spaces\sss $f_v\colon V\toto\cc$\sss 
taking\qss $v$\qss to\qss $1\in\cc$\nnsp.\qss
For a non-zero vectors\sss $v\in V$\sss
we will denote by\qss $q_{\fff v}$\qss the pull-back quadratic form\sss 
$f_v^*\qff \mathbold{n}$\nnsp.\qss
If\dss $V\qff =\qff \cc$\sss and\sss $v\qff =\qff 1$\nnsp,\qss 
then\sss $f_1\qff =\qff f_v\qff =\qff \id_{\dff \cc}$\sss 
and\dss hence\sss $q_{\fff 1}\qff =\qff \mathbold{n}$\nnsp.\oss

\mypar{Lemma.}{change-of-f} \emph{Let\qss $u\fff,\pff v$\qss be two non-zero vectors in\qss $V$\dnsp.\oss
Let $a$ be the unique complex number such that\qss $v\off =\off a\fff u$\nnsp.\qff\oss 
Then\oss 
$\dis
f_u\fff(v)\off =\off a$\oss and\oss 
$\dis
f_u\off =\off a\fff f_v\off =\off m_{\fff a}\circ f_v$\nsp.\oss}

\proof Obviously,\oss 
$f_u\fff (v)\qff =\qff f_u\fff (a\fff u)\qff =\qff a\fff f_u(u)\qff =\qff a\fff 1\qff =\qff a$\oss 
and\oss $a\fff f_v\fff (v)\qff =\qff a\fff 1\qff =\qff a$\nnsp.\oss
Hence $f_u$ and $a\fff f_v$ agree on the basis of $V$ formed by the vector $v$\dnsp.\oss 
It follows that\qss $f_u\qff =\qff a\fff f_v$\nnsp.\qff\oss 
Since\qss $a\fff f_v\qff =\qff m_{\fff a}\circ f_v$\nnsp,\oss
this completes the proof.\oss  \eproof

\mypar{Lemma.}{change-of-vector} \emph{Let\qss $u\fff,\pff v$ be two non-zero vectors in $V$\dnsp,\oss
and let\qss $a$\qss be the unique complex number such that\qss $v\qff =\qff a\fff u$\dnsp.\oss
Then\oss 
$\dis
q_{\fff v}\off =\off |\dff a\dff|^{\dff -\dff 2}\qff q_{\fff u}$\nsp.\oss}

\proof Let us consider the case $V\qff =\qff \cc$\nnsp,\qss 
$u\qff =\qff 1$\qss first\halfff.\oss 
In this case\qss $v\qff =\qff a$
and we need to prove that\qss 
$\dis
q_{\fff a}
\off  =\off  
|\dff a\dff|^{\dff -\dff 2}\dff\mathbold{n}$\nnsp.\qff\oss
Let\qss $b\qff =\qff a^{\dff -\dff 1}$\dnsp.\qff\oss 
Then\qss $m_{\dff b}(a)\qff =\qff 1$\dnsp,\oss 
and hence\oss 
$\dis
f_a\off =\off m_{\dff b}$\dnsp.\qff\oss  
It follows that\oss
$\dis
q_{\fff a}\off =\off m_{\dff b}^*\qff \mathbold{n}$\nnsp.\qff\oss
By Section\qss \ref{quadratic forms}\oss 
$\dis
m_{\dff b}^*\qff \mathbold{n}
\off =\off
|\dff b\dff|^{\fff 2}\dff\mathbold{n}$\nnsp,\oss
and hence\oss\vspace{3pt}
\[
\quad
q_{\fff a}
\off =\off 
m_{\dff b}^*\qff \mathbold{n}
\off =\off
|\dff b\dff|^{\fff 2}\dff\mathbold{n}
\off =\off
|\dff a\dff|^{\fff -\fff 2}\dff\mathbold{n}\dff.
\]

\vspace{-12pt}\vspace{3pt}
This proves the lemma in the case $V\qff =\qff \cc$\nnsp,\qss $u\qff =\qff 1$\dnsp.\oss

In the general case\qss $f_u\fff (v)\qff =\qff a$\qss
by Lemma\qss \ref{change-of-f}\qss and\qss
$f_a\fff (a)\qff =\qff 1$\qss 
by the definition of\qss $f_a$\dnsp.\oss
It follows that\oss 
$\dis
f_a\circ f_u\dff (v)
\off  = \off 
f_a\fff (f_u\fff (v))
\off  =\off  
f_a\fff (a)
\off  =\off 
1$\oss 
and hence\qss  
$f_v\qff =\qff f_a\circ f_u$\qss and\oss\vspace{3pt}
\[
\quad
q_{\fff v} 
\off =\off
f_v^*\qff \mathbold{n} 
\off =\off 
(f_a\circ f_u)^*\qff \mathbold{n} 
\off =\off 
f_u^*\qff (f_a^*\qff \mathbold{n}) 
\off =\off 
f_u^*\qff q_{\fff a}\dff.
\]

\vspace{-12pt}\vspace{3pt}
By the already proved special case,\oss 
$\dis
q_{\fff a}
\off =\off 
|\dff a\dff |^{\fff -\dff 2}\qff q_{\fff 1}
\off =\off
|\dff a\dff |^{\fff -\dff 2}\qff \mathbold{n}$\dnsp.\qff\oss
It follows that\vspace{3pt}
\[
\quad
q_{\fff v} 
\off =\off 
f_u^*\qff q_{\fff a}
\off =\off 
f_u^*\qff (\dff |\dff a\dff |^{\fff -\dff 2}\qff \mathbold{n} \dff) 
\off =\off 
|\dff a\dff |^{\fff -\dff 2}\qff f_u^*\qff \mathbold{n} 
\off =\off 
|\dff a\dff |^{\fff -\dff 2}\qff q_{\fff u}\dff.\pquad\mbox{\eproof}
\]

\vspace{-12pt}\vspace{3pt}
\mypar{Corollary.}{conf-correct} \emph{The conformal class of 
$q_{\fff v}$ does not depend on the choice of non-zero\qss $v\in V$\dnsp.\oss}

\proof Indeed\halfff,\oss if\qss $u\fff,\pff v$\qss are non-zero vectors in $V$\dnsp,\oss
then\qss $v\qff =\qff a\fff u$\qss for some non-zero\qss $a\dff\in\dff\cc$\qss and hence\oss
$\dis
q_{\fff v}
\off =\off 
|\dff a\dff|^{\dff -\dff 2}\qff q_{\fff u}$\nnsp.\oss 
Since $|\dff a\dff|^{\dff -\dff 2}$ is a positive real number\halfff,\oss
the quadratic forms\qss 
$q_{\fff w}$\qss and\qss $q_{\fff v}$\qss
have the same conformal class.\oss  \eproof

\myuppar{The Poincar\'e conformal dilatation of real linear maps.}
Let\sss $V,\off W$\sss be  complex vector spaces of dimension $1$\nnsp.\qss
Let\sss $T\colon V\toto W$\sss be a real-linear map.\oss 
Let us choose non-zero vectors $v\in V$\dnsp,\dss $w\in W$\dnsp.\qss
If we identify both $V$ and $W$ with $\cc$ by the maps\qss 
$f_v\colon V\to\cc$\qss and\qss $f_w\colon W\to\cc$\qss respectively,\qss 
the map $T$ will turn into a real linear map $\cc\toto\cc$\nnsp,\qss
allowing to speak about its Poincar\'e complex dilatation 
as defined in Section\qss \ref{linear maps}.\oss 
Using these identification is equivalent to replacing $T$ by
the map\sss $f_w\circ T\circ f_v^{\dff -\dff 1}$\sss considered as a map\sss $\cc\toto\cc$\dnsp.\qss
Naively,\oss one may try to define the Poincar\'e conformal dilatation of $T$ 
as the Poincar\'e complex dilatation $\pi_{\dff T}\dff(v\fff,\pff w)$ 
of\dss $f_w\circ T\circ f_v^{\dff -\dff 1} $\dnsp.\qss
As we will\sss see in a moment,\dss
$\pi_{\dff T}\dff(v\fff,\pff w)$
does not\sss depends on $w$\nnsp,\qss
but\qss \emph{does depends on}\dss $v$\nnsp.\oss

\mypar{Lemma.}{no-target} \emph{\dnsp$\pi_{\dff T}\dff(v\fff,\pff w)$\pss 
does not depend on\qss $w$\sss
and\dss hence can be denoted\dss by $\pi_{\dff T}\dff(v)$\nnsp.\qss}

\proof It is sufficient to prove that the conformal class of the form\qss
$(f_w\circ T\circ f_v^{\dff -\dff 1})^*\qff n$\qss does not depend on $w$\dnsp.\qff\oss
Since\sss
$(\fff f_w\circ T\circ f_v^{\dff -\dff 1})^*\qff n 
\off =\off 
(T\circ f_v^{\dff -\dff 1})^*\qff(\fff f_w^*\qff n)
\off =\off
(T\circ f_v^{\dff -\dff 1})^*\qff q_{\fff w}$\nsp,\qss
Corollary\qss \ref{conf-correct}\qss implies that this is indeed the case.\oss  \eproof

\mypar{Lemma.}{source-change} \emph{Let\dss $u\fff,\pff v$ be two non-zero vectors in $V$\dnsp,\oss
and let\dss $a$\sss be the unique complex number such that\sss $v\qff =\qff a\fff u$\nsp.\qff\oss
Then\dss
$\pi_{\dff T}\dff(v)\off =\off (\overline{a}\fff/a)\qff \pi_{\dff T}\dff(u)$\dnsp.\oss}

\proof Let\qss $w\in W$\dnsp,\dss $w\qff \neq\qff 0$\nnsp.\qss 
By the definition,\dss 
$\pi_{\dff T}\dff(v)$\sss and\sss $\pi_{\dff T}\dff(u)$\sss are equal to
the Poincar\'e complex dilatations of the maps\sss
$f_w\circ T\circ f_v^{\dff -\dff 1}$\sss and\sss
$f_w\circ T\circ f_u^{\dff -\dff 1}$\sss respectively.\qss
Lemma\qss \ref{change-of-f}\qss implies\sss that\sss
$f_u
\off =\off 
m_{\fff a}\circ f_v$\oss
and\dss hence\sss 
$f_v^{\dff -\dff 1}
\off =\off 
f_u^{\dff -\dff 1}\circ m_{\fff a}$\oss and\vspace{4.5pt}
\[
\quad
f_w\circ T\circ f_v^{\dff -\dff 1} 
\off =\off
f_w\circ T\circ f_u^{\dff -\dff 1}\circ m_{\fff a}
\off =\off
\left(\fff f_w\circ T\circ f_u^{\dff -\dff 1}\fff\right)\circ m_{\fff a}\dff.
\]

\vspace{-12pt}\vspace{4.5pt}
It follows that\vspace{4.5pt}
\begin{equation}
\label{source-compare-maps}
\quad
(f_w\circ T\circ f_v^{\dff -\dff 1})^*\qff \mathbold{n} 
\off =\off
((f_w\circ T\circ f_u^{\dff -\dff 1})\circ m_{\dff a})^*\qff \mathbold{n} 
\off =\off 
m_{\dff a}^*\qff (\dff (f_w\circ T\circ f_u^{\dff -\dff 1})^*\qff \mathbold{n}\dff)\fff.
\end{equation}

\vspace{-12pt}\vspace{4.5pt}
Let $q\qff =\qff (f_w\circ T\circ f_u^{\dff -\dff 1})^*\qff \mathbold{n}$\nnsp.\oss
Then\qss (\ref{source-compare-maps})\qss implies that\qss
$\dis
(f_w\circ T\circ f_u^{\dff -\dff 1})^*\qff \mathbold{n}
\qff =\qff
m_{\dff a}^*\qff q$\nnsp,\oss 
and\dss hence\vspace{4.5pt}
\[
\quad
\pi_{\dff T}\dff(v)
\off =\off 
P\fff(\fff m_{\dff a}^*\qff q\fff )
\off =\off
(\overline{a}\fff/a)\dff P\fff(\fff q\fff )
\off =\off
(\overline{a}\fff/a)\dff \pi_{\dff T}\dff(u)\dff,
\] 

\vspace{-12pt}\vspace{4.5pt}
where\sss the second equality\sss follows from\qss (\ref{p-rotations}).\oss  \eproof

\myuppar{The dual and the conjugate vector spaces.}
Lemma\qss \ref{source-change}\qss means\sss that\sss the Poincar\'e complex dilatation\sss
$\pi_{\dff T}\dff(v)$\sss is\sss a quantity depending on\sss the choice of\dss
a basis of\sss $V$\dnsp.\oss
Moreover\halfff,\qss this quantity depends on\sss the choice in a way\sss
making\sss it\sss a\qss \emph{tensor}\qss in a classical\sss sense.\oss
Now we will\sss explain\sss what\sss this means in\sss the modern\sss language.\qss

Let $V$ be a complex vector space.\oss 
Recall\dss that\sss its\qss \emph{dual}\qss vector space $V^*$
has as vectors the complex linear maps $V\toto\cc$\nnsp.\oss
If\dss $a\dff\in\dff \cc$\sss and\sss $f\dff\in\dff V^*$\dnsp,\qss
then\sss $a\fff f$\sss is defined by\sss 
$\dis
a\dff f\fff(v)
\off =\off 
a\dff(f\fff(v))
\off =\off
f\fff(a\fff v)$\nnsp.\oss
The\qss \emph{conjugate}\qss vector space $\overline{V}$ of $V$
is equal to $V$ as a real vector space,\oss
but the multiplication by $a\dff\in\dff\cc$ in $\overline{V}$ is defined 
as the multiplication by $\overline{a}$ in $V$\dnsp.\qss
We will denote\sss this multiplication by
$v\qff\longmapsto\qff a\cdot v$\dnsp,\qss 
so that $a\cdot v\qff =\qff \overline{a}\dff v$\dnsp.\oss 

Suppose now that the complex dimension of $V$ is $1$\dnsp.\oss 
Every non-zero\qss $v\dff\in\dff V$\qss forms
a basis of\dss both $V$ and $\overline{V}$\dnsp.\oss
The dual basis of\sss $\overline{V}^{\qff *}$ 
consists of the map $v^{\dff\bullet}\colon \overline{V}\toto \cc$ 
defined by $v^{\dff\bullet}(w)\qff =\qff \overline{f_v\fff(w)}$\dnsp.\qss
Indeed,\dss
$v^{\dff\bullet}(v)\qff =\qff \overline{f_v\fff(v)}\qff =\qff 1$\nnsp,\qss
and\vspace{4.5pt}
\[
\quad
v^{\dff\bullet}(a\cdot w)
\off =\off
\overline{f_v\fff(a\cdot w)}
\off =\off
\overline{f_v\fff(\overline{a}\dff w)}
\off =\off
\overline{\overline{a}\qff f_v\fff(w)}
\off =\off
a\dff \overline{f_v\fff(w)}
\off =\off
a\dff v^{\dff\bullet}(w)
\dff,
\]

\vspace{-12pt}\vspace{4.5pt}
i.e.\dss $v^{\dff\bullet}$\sss is\sss complex-linear\halfff.\qss
The tensor product $\overline{V}^{\qff *}\otimes V$ is also of complex dimension $1$\dnsp.\oss
Any non-zero vector $v\dff\in\dff V$ leads to a basis of $\overline{V}^{\qff *}\otimes V$ 
consisting of the vector $v^{\dff\bullet}\otimes v$\nnsp.\oss

\mypar{Theorem.}{basis-independence} \emph{Let\qss $V,\off W$\qss 
be complex vector spaces of dimension\dss $1$\dnsp.\oss
Let\qss $T\colon V\toto W$\qss be a real linear map.\oss
Then\oss
$\dis
\pi_{\dff T}\fff(v) \left(\fff v^{\fff\bullet}\otimes v\fff\right)$\oss
does not depend on the choice of a non-zero\qss $v\in V$\dnsp.\oss}

\proof Suppose that\qss $v\fff,\pff u\dff\in\dff V$\qss are non-zero.\oss
Then\qss $v\qff =\qff a\fff u$\qss for some\qss $a\dff\in\dff \cc^*$\dnsp.\oss
By Lemma\qss \ref{change-of-f}\qss $f_u\qff =\qff a\dff f_v$\nsp.\oss
It follows that\oss 
$\dis
u^{\dff\bullet}
\off =\off 
\overline{a}\qff v^{\dff\bullet}$\nsp.\qff\oss
Therefore\oss
$\dis
\left(\dff\overline{a}\dff\right)^{\fff -\fff 1} u^{\dff\bullet}
\off =\off 
v^{\dff\bullet}$\oss
and
\[
\quad
v^{\dff\bullet}\otimes v
\off =\off
\left(\dff(\overline{a})^{\fff -\fff 1} u^{\dff\bullet}\dff\right)\otimes {a\fff u}
\off =\off
\left( a/\dff\overline{a}\fff\right) \left(\fff u^{\dff\bullet}\otimes u\fff\right)
\]
By Lemma\qss \ref{source-change}\oss
$\dis
\pi_{\dff T}\dff(v)\off =\off \left(\fff\overline{a}/a\right)\qff \pi_{\dff T}\dff (u)$\dnsp.\oss
It follows that\vspace*{\smallskipamount}
\[
\quad
\pi_{\dff T}\fff(v) \left(\fff v^{\dff\bullet}\otimes v\fff\right)
\off =\off
(\fff\overline{a}/a)\qff \pi_{\dff T}\dff(u) \left(\fff v^{\dff\bullet}\otimes v\fff\right)
\]

\vspace*{-3\bigskipamount}
\[
\quad
\phantom{\pi_{\dff T}\fff(v) \left(\dff v^{\dff\bullet}\otimes v\fff\right)
\off }
=\off
(\fff\overline{a}/a)\qff\pi_{\dff T}\dff(u) 
\left( a/\dff\overline{a}\fff\right) \left(\fff u^{\dff\bullet}\otimes u\fff\right)
\off =\off
\pi_{\dff T}\fff(u) \left(\dff u^{\dff\bullet}\otimes u\fff\right).
\]

\vspace*{-0.75\bigskipamount}
Therefore,\oss 
$\dis
\pi_{\dff T}\fff(v) \left(\fff v^{\fff\bullet}\otimes v\fff\right)$\qss 
indeed does not depend on the choice of $v$\nnsp.\oss  \eproof

\myuppar{Poincar\'e\dss conformal\dss dilatation as a tensor\halfff.}
In view of Theorem\qss \ref{basis-independence}\qss we can define
the\qss \emph{Poincar\'e\dss conformal\dss dilatation}\qss 
$\pi_{\dff T}$\qss of\dss $T$\dss by the formula\sss
$\pi_{\dff T}
\off =\off
\pi_{\dff T}\fff(v) \left(\fff v^{\fff\bullet}\otimes v\fff\right)
\qff \in\qff
\overline{V}^{\qff *}\otimes V$\sss
for any non-zero\sss $v\dff\in\dff V$\dnsp.\qss
Suppose that\sss $V\qff =\qff W\qff =\qff \cc$\nnsp.\qss
In this case the identification of\dss 
$\overline{\cc}^{\qff *}\otimes \cc$\sss with\sss $\ccc$\sss
by the isomorphism taking $1^{\dff\bullet}\otimes 1$ to $1$\sss
turns\sss
$\pi_{\dff T}\dff\in\dff \overline{\cc}^{\qff *}\otimes \cc$\sss
into the Poincar\'e conformal dilations\sss $\pi_{\dff T}$\sss
from Section\qss \ref{linear maps},\qss as one easily checks.\qss 
This easily\sss implies\sss that\sss the Poincar\'e conformal dilation 
agrees with\sss the classical\sss one in\sss general.\qss

\begin{flushright}

January\qss 21,\oss 2017.\oss September\qss 10,\oss 2023\qss --\qss The present\sss version.
 
https\halfff:/\!/{\hnsp}nikolaivivanov.com

E-mail\halfff:\oss nikolai.v.ivanov{\fff}@{\fff}icloud.com,\oss ivanov{\fff}@{\fff}msu.edu

\end{flushright}

\end{document}